\documentclass[a4paper, 12pt]{article}
\usepackage{amsmath}
\usepackage{amsfonts}
\usepackage{amssymb}
\usepackage{amscd}
\usepackage{amsthm}

\renewcommand{\P}{{\mathbb P}}
\newcommand{\Z}{{\mathbb Z}}
\newcommand{\Q}{{\mathbb Q}}
\newcommand{\Ge}{{\geqslant}}

\newcommand{\Ql}{{\mathbb Q}_l}

\newcommand{\R}{{\mathbb R}}
\newcommand{\C}{{\mathbb C}}
\newcommand{\Cb}{{\mathbb C}}

\newcommand{\F}{{\mathcal F}}
\newcommand{\Fb}{{\mathbb F}}

\newcommand{\Dc}{{\mathcal D}}

\newcommand{\Nm}{{\rm Nm}}

\newcommand{\Sc}{{\mathcal S}}

\newcommand{\A}{{\mathbb A}}

\newcommand{\OO}{{\mathcal O}}

\newcommand{\Id}{{\rm Id}}
\newcommand{\Hom}{{\rm Hom}}

\newcommand{\vol}{{\rm vol}}
\newcommand{\Gal}{{\rm Gal}}

\newcommand{\Pic}{{\rm Pic}}

\newcommand{\Supp}{{\rm Supp}}

\newcommand{\Ker}{{\rm Ker}}

\newcommand{\Tr}{{\rm Tr}}
\newcommand{\Fo}{{\rm F}}

\newcommand{\Spec}{{\rm Spec}}

\newcommand{\res}{{\rm res}}

\renewcommand{\dim}{{\rm dim}}

\newcommand{\Div}{{\rm Div}}

\theoremstyle{plain}

\newtheorem{defin-prop}[theor]{Definition-Proposition}

\theoremstyle{remark}

\theoremstyle{definition}

\title{Notes on the Poisson formula}
\author{A.\,N.\,Parshin\footnote{Steklov Mathematical Institute, Russian Academy of Sciences, Gubkina str 8, Moscow 119991, Russia; parshin@mi.ras.ru.}}
\date{}

\begin{document}
\maketitle

\begin{itemize}
\item[\sf]Introduction
\item[\sf]The Tate--Iwasawa method for algebraic curves
\begin{itemize}
\item[\sf]      Preliminaries  and notations
\item[\sf]      Duality and the Fourier transform
\item[\sf]      The Weil problem
     \item[\sf] Regularization
     \item[\sf] Applications
\end{itemize}
\item[\sf]Discrete versions and holomorphic tori
\begin{itemize}
       \item[\sf]Analysis on discrete groups
       \item[\sf]Extended group of divisors
      \item[\sf] Holomorphic tori
      \item[\sf] The Poisson formula and residues
      \item[\sf] Explicit formulas
\end{itemize}
\item[\sf]Trace formula and the Artin representation
\begin{itemize}
      \item[\sf]  Zeta-functions and \'etale cohomology
       \item[\sf] Local Artin representation
       \item[\sf] Relation with the explicit formulas
\end{itemize}
\item[\sf]Number fields (from  $\C^*$  to $\C$)
\item[\sf]Appendix: Irina Rezvjakova.
Residues of the Dedekind zeta functions of number fields
\item[\sf]References
\end{itemize}

These notes are a part of my lectures on representations of adelic groups attached
to two-dimensional schemes.  The main part of the lectures was  focused on the representation theory
for discrete Heisenberg groups and applications to the groups of this type arising from the algebraic surfaces defined
over a finite field. The arithmetic surfaces were not considered except for the  remark that the theory for geometric surfaces
can be transposed to this case without changes if we restrict ourselves to  the scheme part of the arithmetic surface.

 The classical  representation
theory for  discrete groups is not well developed  since these groups are mostly  not of type I. This implies a violation
of the main principles  of  representation theory on  Hilbert spaces: non-uniqueness of  the decomposition into irreducible components;
 bad topology of the unitary dual space; non-existence of characters... (see, for example, \cite{Dix}). According to  the Thoma theorem a discrete group is of type I if and only if  it has an abelian subgroup of finite index.
On the other hand, there exists a theory of smooth representations for  $p$-adic algebraic groups. This theory is also valid for a more general class
of totally  disconnected locally compact groups. The discrete groups are a simple particular case of this class of groups and the general theory delivers a reasonable class of  representations, which are   representations on a vector space without any topology.
The new viewpoint consists in a systematic consideration of  purely algebraic induced  representations in place of   unitary representations on Hilbert spaces.

In our situation, the change in  the class
of representations will cause  the moduli spaces of induced representations  to be  complex-analytic manifolds. Characters do exist and will then be modular forms (see a brief exposition in \cite{P3}). It seems that the more general holomorphic dual
space is more adequate for this class of groups then the standard unitary dual going back to the Pontrjagin duality for abelian groups.

This approach can change even the theory of representations for abelian discrete groups and its application to the arithmetic of
algebraic curves. This point, we hope, will be clear from these notes.

The lectures were given in both divisions of the Steklov Mathematical Institute (Saint-Petersburg, November 17-20,
and Moscow, December 21-30, 2009)\footnote{A substantial part of these notes were also presented in my lectures  in the Humboldt University, Berlin, October 2010.}. I am very thankful to audience, especially to Volodya Zhgun,
Maxim Korolev and Irina Rezvjakova, for   fruitful discussions and to Lawrence Breen who  read the final text  and made many important remarks.
Also, I am grateful to Sergey Gorchinsky and Igor Zhukov and his students for their
help in  preparation of these notes.

\section{Introduction}

The simplest case of the Poisson summation formula is the following one. Let $\R$ be the additive group of
real numbers and $\Z$ be the discrete subgroup of integer  numbers. Denote by $\Sc(\R)$ the space of
smooth  functions on  $\R$, rapidly  decreasing together  with all its derivatives (=  Schwartz functions \cite{S}).
We have the Fourier transform
$$
\Fo:  \Sc(\R) \rightarrow \Sc(\R),
$$
where for $f  \in \Sc(\R)$ and $y \in \R$
$$
\Fo(f)(y)  = \int \limits_{\R}\exp (2\pi ixy)f(x)dx.
$$
The Poisson formula reads (see \cite{GSh})
$$
\sum_{n \in \Z} f(n) =  \sum_{n \in \Z} (\Fo f)(n)
$$
The main ingredients of this result are the following items:
\begin{itemize}
\item a  locally compact commutative topological group $G$
\item  its Pontryagin dual group $\check G$
\item  a discrete subgroup $\Gamma \subset G$
\item  a space of functions $\F(G)$
\item  a Haar measure  $\mu \in \mu(G)$ on   $G$ such that $\mu(G/\Gamma) < \infty$
\item  the Fourier transform
$$
\Fo:  \F(G) \rightarrow \F(\check G).
$$
\end{itemize}
In the  simplest case $ G = \R = \check G,~ \Gamma = \Z, ~\mu = dx$ and $\F(G) = \Sc(\R)$.
In the  more general situation, one has
\begin{equation}\label{eq-pf}
\sum_{x \in \Gamma} f(x) =  \mu(G/\Gamma)^{-1}\sum_{y \in \Gamma^{\perp}} (\Fo f)(y),
\end{equation}
where $\Gamma^{\perp} = \{ y \in G : \langle y, x\rangle  = 1~ \mbox{for all } x \in \Gamma\}$ is the annihilator (or orthogonal) of
$\Gamma$,
$$
\Fo (f)(y)  = \int \limits_{G}\langle y, x\rangle f(x)d\mu(x),
$$
 $\langle -, -\rangle: G \times \check G \rightarrow {\bf U}(1)$\footnote{${\bf U}(1)$ is the group of complex numbers with modulus equals to 1.} is the canonical pairing of $G$ with $\check G$ and   $\F(G)$ is an appropriate
 space of functions.

Let us outline a proof of the Poisson formula. Take the function $\varphi (g): = \sum\nolimits_{\gamma \in \Gamma} f(g\gamma), ~g \in G/\Gamma$ and integrate it over
$G/\Gamma$ with characters $\chi \in \check G$  trivial on $\Gamma$:
$$
(\Fo\varphi) (\chi) = \int\limits_{G/\Gamma}\chi(g)\varphi(g)d\mu(g).
$$
The measure on $G/\Gamma$ is basically the same as the original measure on $G$ and the set of characters
will be exactly $\Gamma^{\perp} = \mbox{dual to}~G/\Gamma$ . Thus, the integral is the Fourier transform on the group
$G/\Gamma$. Assuming that the integral (= sum) of this function over $\Gamma^{\perp}$  converges,  we have
$$
\sum_{\chi \in \Gamma^{\perp}}  (\Fo\varphi) (\chi)  = \varphi (e) = \sum_{\gamma \in \Gamma}f(\gamma),
$$
where the first equality follows from the inversion formula $\Fo \circ \Fo  = (-1)^*$ for the Fourier transform on the group
$G/\Gamma$ and it's dual group. On the other hand,
interchanging in the first integral summation (defining the $\varphi$) and integration, we have
$$
 (\Fo\varphi) (\chi) =  \sum_{\gamma \in \Gamma}\int_{G/\Gamma}f(g\gamma)\chi(g)d\mu(g) =
 \sum_{\gamma \in \Gamma}\int_{G/\Gamma}f(g\gamma)\chi(g\gamma)d\mu(g) =
 \int_{G}f(g)\chi(g)d\mu(g) =  (\Fo f)(\chi).
 $$
Combining all this, we get the Poisson formula.

Beyond the equality (\ref{eq-pf}) and the proof, one can find some functorialities  which will be valid in a more
general situation, when one has an arbitrary closed subgroup $H \subset G$ instead of the discrete subgroup
$\Gamma$. We will denote by $\mu(G)$ the space of Haar measures on a group $G$. Thus, $\mu(G) \cong \C$
and the Fourier transform is a canonical map
$$
\Fo_G: \F(G)\otimes \mu(G) \rightarrow \F(\check G)
$$
defined for any locally compact group $G$.

If our group $G$ is totally disconnected then we set
$$
\F(G) = {\cal D}(G) =  \{\mbox{locally constant functions with compact support}\}.
$$
These are the classical spaces introduced by F. Bruhat \cite{Br}.
For connected Lie groups we can use the space  $\Sc(G)$ of  the Schwartz functions \cite{S}
and its dual $\Sc'(G)$. The general case can be easily reduced to these two (see \cite{Br}).
The Fourier transform for the distribution spaces $\F'(G)$ such as  ${\cal D}'(G) =
\{\mbox{dual to}~{\cal D}(G) \mbox{, i.e. all distributions}\}$
and  $\Sc'(G)$ can be defined by duality and will be a map
$\Fo: \F'(G)\otimes \mu(G)^{-1} \rightarrow \F'(\check G)$, where $\mu(G)^{-1}$
is the dual 1-dimensional space. The harmonic analysis includes definitions of direct and inverse
images for the functions and the distributions in a category of  locally compact commutative
groups. The Fourier transform will interchange the direct and inverse images.

Suppose  we are given a closed subgroup $i \, : \, H \rightarrow G$. Denote by $\pi$ the canonical
 map from $G$ onto $G/H$, by $\alpha$ the embedding of trivial group $\{ e\}$ into  $G/H$
 and by $\beta$ the surjection of $H^{\perp}$ onto $\{ e\}$. Here,  ${H^{\perp} \subset \check{G}}$
 is the annihilator of $H$ in $\check G$ defined as above.

Then the diagram
\begin{equation}\label{d}
\begin{CD}
{\cal F}(G)\otimes \mu(G)   @> \pi_* >>     {\cal F}(G/H)\otimes \mu(G/H)  @> \alpha^*\otimes \Id >> \mu(G/H) \\
@ V \Fo_G VV          @V \Fo_{G/H} VV    @V \cong VV \\
{\cal F}(\check G)     @> \check \pi^* >>      {\cal F}(H^{\perp})  @> \beta_* >> \mu(H^{\perp})^{-1}
\end{CD}
\end{equation}
commutes. Here the direct image $ \pi_*$ is  the  integral over the fibers of $\pi$ and the inverse image
$ \pi^*$ is restriction to the subgroup $H^{\perp}$.  To define the integral one must use the canonical decomposition
$\mu(G) = \mu(H)\otimes\mu(G/H)$. Since the Haar measure is translation
invariant, this gives measures not only on $H = \Ker \pi$ but on all fibers of $\pi$.
The map $\alpha^*$ is the evaluation at the point $e$ and
$ \beta_*$ is the integral over $H^{\perp}$. Finally, the canonical isomorphism between measure spaces follows from the equalities
$\check{(G/H)}  = H^{\perp}$ and $\mu(\mbox{a group})\otimes \mu(\mbox{its dual group}) = \C$.

Then we have the  following Poisson formula
\begin{equation}\label{eq-pf2}
\Fo(\delta_{H,~\mu_0} \otimes \mu^{-1}) = \delta_{{H^{\perp}},~ \mu^{-1}/\mu^{-1}_0}
\end{equation}
for any closed subgroup $i \, : \, H \rightarrow G$.
Here
\begin{itemize}
\item   $\mu_0 \in \mu(H) \subset {\cal D}'(H)$
\item $\mu \in \mu(G)$
\item   $\delta_{H,\mu_0} = i_*(\mu_0),\quad i_*: {\cal D}'(H) \rightarrow {\cal D}'(G)$
\end{itemize}

We use that, by definition of measures, $\mu(H^{\perp}) = \mu(G/H)^{-1} =   \mu(G)^{-1}\mu(H)$.

If the subgroup $H$ is discrete then  one can take for  $\mu_0 \in \mu(H)$ the canonical
measure on  $H$,  and the $\delta_{H,\mu_0}$ will be exactly the sum over   $H$ and we arrive at  the formula (\ref{eq-pf}).
Concerning the general formalism of harmonic analysis over $n$-dimensional local fields and adelic groups ($n = 0, 1, 2$),  see
\cite{P1, OP1, OP2}. The basic problem solved there was to extend the classical analysis known for locally compact (and first of all for finite)
groups to the case of   fields and groups arising from two-dimensional schemes. These fields and groups  are no longer  locally compact.

\bigskip

We will now briefly outline one of the main applications of the Poisson formula to investigation
of the analytic structure of Riemann's zeta function $\zeta(s)$. We know that
$$
 \zeta(s) = \sum_{n \ge 1}n^{-s}  = \prod_p (1 - p^{-s})^{-1} = \prod_p\zeta_p(s),
$$
where $p$ runs through the prime numbers, $s \in \C$ and both expressions are convergent for
${\rm Re} s > 1$. The sum is the Dirichlet series and the product is the Euler product.
According to the general principle of arithmetic (see, for example,  \cite{P2}) one has to extend the definition of $\zeta(s)$ by
a factor $\zeta_{\infty}(s)$ at infinity. The new $\zeta$-function will no longer correspond to the affine scheme
$\Spec (\Z)$ but to a "complete" scheme  $C = \Spec (\Z)  \cup \infty$. Just as the primes $p$ correspond to
$p$-adic completions $\Q_p$ of the field  $\Q$ ~(= field of rational functions on $\Spec (\Z)$),  the
$\infty$ corresponds to the embedding of  $\Q$ into $\R$. The missing factor was introduced many
years before the ideas mentioned here  appeared. Namely, we have
$$
\zeta_{\infty}(s)  =  \pi^{-s/2}\Gamma(s/2)
$$
and
$$
 \zeta_C(s) =  \prod_{x \in C}\zeta_x(s).
$$
Now, let $t \in \R^*$,  $ f_t(x): = \exp(-\pi x^2t) \in  \Sc(\R)$
and $f(t) = \sum_{n \geq 1} \exp (-\pi n^2t) = \frac{1}{2}(\sum_{n \in Z}f_t(n) - 1)$.
Then
$$
F(f_t) = t^{-1/2}f_{t^{-1}}.
$$
For the gamma function, we know that
$\Gamma(s) = \int \limits_{0}^{\infty} \exp (-t)t^s \frac{dt}{t}$.
Then, if we use the Dirichlet series for $\zeta(s)$ and make  change of variables
$t \mapsto \pi n^2t$, we   easily get that
$$
\zeta_C(s) = \int \limits_{\R^*_+}t^{s/2} f(t) \frac{dt}{t}.
$$
Breaking the domain of integration $\R^*_+$ into two pieces, $t < 1$ and $t > 1$ and applying the simplest Poisson formula
to the function $ f_t(x)$, we arrive at  the following
expression for the zeta-function\footnote{We recommend  that the reader  make all the computations by himself and then  compares
them with the computations for  zeta-functions of curves in the next section. A completely unified  exposition of both cases, for number fields and
for curves, is given in Weil's book \cite{W2}[ch. VII, \S 5]}:
\begin{equation}\label{eq-rz}
\zeta_C(s) =  -  \frac{1}{s} +  \frac{1}{s-1}  + \int \limits_1^{\infty} f(t)[ t^{\frac{s}{2}} + t^{\frac{1-s}{2}}]\frac{dt}{t},
\end{equation}
where the integral over $t > 1$ converges for all $s \in \C$. This gives an analytic continuation of $\zeta_C(s)$ and $\zeta(s)$ to the whole $s$-plane $\C$ and shows that the singularities
are  poles at $s = 0$ and $s = 1$.  The functional equation has the following form
$$
\zeta_C(1 - s) = \zeta_C( s).
$$
This approach belongs to Riemann \cite{R} and  was extended to zeta- and $L$-functions of the fields of algebraic numbers (= finite extensions of
$\Q$) by Hecke \cite{H1, H2}, who  used a functional equation for the theta-functions (see the appendix below). In our case, the function $f(t)$ is a simple example of a theta-function. In the next sections we will study this construction for the parallel case of algebraic curves defined over a finite field. After
that we will return  to the number fields.

\section{The Tate--Iwasawa method for algebraic curves}

\subsection{Preliminaries and notations}

Let $C$ be a smooth projective curve over a finite field $k=\Fb_q$.
Let $K=k(C)$ be the field of rational functions on $C$. Denote by
$\A$ the ring of adeles on $C$, that is
$$
\A:=\prod_{x\in C}{'K_x},
$$
where $x$ runs through all closed points of $C$, $K_x$ is the
completion of the field $K$ at a point $x$, and the restricted (adelic)
product is taken with respect to the discrete valuation subrings
$\hat\OO_x$ in $K_x$. A local parameter $t_x$ at a point $x\in C$
determines isomorphisms
$$
K_x\cong k(x)((t_x)),\quad \hat\OO_x\cong k(x)[[t_x]],
$$
where $k(x)$ is the residue field of $x$. Furthermore, the field $K_x$
is locally compact with respect to the topology defined by the discrete
valuation $\nu_x$ and $\hat\OO_x$ is an open compact subring in $K_x$. It
follows that $\A$ is a topological ring and is locally compact.
Recall that the diagonal embedding
$K\hookrightarrow \A$ induces the discrete topology on $K$, while the
quotient $\A/K$ is compact. Set
$\OO:=\prod\limits_{x\in C}\OO_x$ and for a divisor $D=\sum\limits_{x\in C} n_x\cdot x = \sum\limits_{x\in C} \nu_x(D)\cdot x$ on
$C$, set
\begin{eqnarray}
K_x(D):=&\{a_x \in K_x,|\, \nu_x(a_x)+n_x\ge 0\},\\
\A(D):=&\{(a_x)\in\A\,|\, a_x \in K_x(D),~ \mbox{for all}~ x~ \in C\}.\label{eq-A(D)}
\end{eqnarray}
In particular, $\OO=\A(0)$. Note that $\A(D)$ is a compact subset of $\A$.

Let $\Dc(\A)$ be the space of complex-valued
Bruhat--Schwartz functions on $\A$, that is, complex valued locally
constant functions with compact support. We consider  $\Dc(\A)$  just as a  vector space, without any topology.
For each divisor $D$ on $C$, the delta-function $\delta_{\A(D)}$
of $\A(D)$ belongs to $\Dc(\A)$. By definition,
$\delta_{\A(D)}(a)=1$ if $a\in \A(D)$ and $\delta_{\A(D)}(a)=0$ otherwise.

An important fact is that the space $\Dc(\A)$ is generated by finite linear combinations of the functions
$\delta_{\A(D)}$\footnote{This fact is well known \cite{W2}.  It follows, for example, from the  presentation of the space
$\Dc(\A)$ as a double inductive limit of the spaces ${\cal F}(\A(D)/\A(D'))$ of all functions on finite sets  $\A(D)/\A(D')$ \cite{P1, OP1}.}.

Let $\Dc'(\A)$ be the
vector space of all complex valued  linear functionals on $\Dc(\A)$, that is, the space of all distributions. An example of a distribution from $\Dc'(\A)$ is the delta-function $\delta_0$, $\delta_0(f)=f(0)$ for any function $f\in\Dc(\A)$. Another example is the delta-function $\delta_K$ of the discrete subgroup $K\subset \A$: for any function $f\in\Dc(\A)$, we have
$$
\delta_K(f)=\sum_{a\in K}f(a).
$$

Denote by $\A^*$ the group of ideles on $C$, that is, invertible
elements in the ring $\A$. We have
$$
\A^*=\prod_{x\in C}{'K_x^*},
$$
where $K_x^*$ is the multiplicative group of the field $K_x$,
the adelic  product is taken with respect to the subgroups $\hat\OO_x^*$ in $K_x^*$, and
$\hat\OO^*_x$ is the group of invertible elements in the ring
$\hat\OO_x$. As above, $K_x^*$ is a locally
compact topological group, $\hat\OO_x^*$ is an open compact subgroup
in $K_x^*$, and $\A^*$ is a locally compact topological group. It can be directly shown that the
embedding $\A^*\subset \A$ has a dense image. Furthermore, there is a
well-defined degree homomorphism
$$
\deg:\A^*\to \Z,\quad
(a_x)\mapsto\deg(a_x):=\sum_{x\in C}\nu_x(a_x)\deg(x),
$$
where $\deg(x)$ is the  degree of the residue field $k(x)$ over $k=\Fb_q$ for any closed point $x\in C$. We have $q^x: = \# k(x) = q^{\deg(x)}$.
Note that the map $\deg$ is surjective\footnote{This is not a trival claim. It means that any curve $C$ contains  a divisor of degree 1! but not necessarily a point
$x$ of degree 1, i.e. rational over $\Fb_q$. Indeed,  the variety $\Pic^1(C)$ parametrize  classes of divisors of degree 1. The variety will be
a torseur for the algebraic group $\Pic^0(C)$ and any torseur for a group over a finite field has a rational point \cite{Se1}[ch. VI, \S 1, n 4]}.
Let $(\A^*)^{(1)}$ be the kernel of the degree homomorphism. Recall
that the diagonal embedding $K^*\hookrightarrow \A^*$ induces a
discrete topology on $K^*$ and  the restriction of the degree to
$K^*$ is identically zero. An important result is that the quotient $(\A^*)^{(1)}/K^*$ is
compact (this corresponds to the fact that degree zero line bundles on $C$ are parameterized by the $\Fb_q$-points of the Jacobian variety, which is a finite set).

Note that each idele $g\in\A^*$ defines distributions $\delta_{gK} \in \Dc'(\A)$ and $\delta_{gK^*} \in \Dc'(\A)$: for any function $f\in\Dc(\A)$, we have
$$
\delta_{gK}(f)=\sum_{a\in gK}f(a),\quad \delta_{gK^*}(f)=\sum_{a\in gK^*}f(a)=\delta_{gK}(f)-\delta_0(f).
$$

We will consider continuous complex-valued characters of
the group $\A^*$ whose restriction to $K^*$ is trivial:
$$
\chi :\A^* \longrightarrow \Cb^*, ~ \chi|_{K^*}\equiv 1.
$$
In particular, the norm homomorphism $|\cdot|:=q^{-\deg(\cdot)}$ is a character of this type.
We will use the following involution on characters:
$$
\breve{\chi}:=|\cdot|\chi^{-1}.
$$

Since the extension of topological commutative groups
$$
1\to (\A^*)^{(1)}/K^*\longrightarrow \A^*/K^*\stackrel{\deg}\longrightarrow \Z\to 0
$$
is  split
 and the group $(\A^*)^{(1)}/K^*$ is compact, any character $\chi$ as above can be written as
$$
\chi=\chi_0 |\cdot|^s
$$
for some $s\in \C$, where $|\chi_0(a)|=1$ for any element $a\in\A^*$. Therefore,
$$
\Hom(\A^*)/K^*, \C^*) = \C^* (= \C/\frac{2\pi i}{\ln q}\Z)\times (\mbox{a discrete group}).
$$
It follows that the real part ${\rm Re}(s)$  of $s$ does not depend on the decomposition above. Denote by ${\rm Re}(\chi)$ the real number ${\rm Re}(s)$.

\subsection{Duality and the Fourier transform}

Let us choose a non-zero rational differential form $\omega$ on $C$ and consider the pairing on $\A$ given by the formula
$$
\langle\cdot,\cdot\rangle:\A\times\A\to {\bf U}(1)
$$
$$
\langle(a_x),(b_x)\rangle:=\sum_{x\in C}\Psi\left(\Tr_{k(x)/k}\,\res_x(a_xb_x\omega)\right),
$$
where ${\bf U}(1)$ can be viewed as the unit circle in $\C$ and $\Psi:k=\Fb_q\to \C^*$ is an injective group homomorphism determined  by a choice of a $q$-th root of unity in $\C$. The pairing above is continuous and non-degenerate, and the locally compact commutative group $\A$ is its own Pontryagin  dual. It is easy to see that the orthogonal of $\A(D)$ with respect to the pairing $\langle\cdot,\cdot\rangle$ is $\A((\omega)-D)$, where $(\omega)$ is the divisor of the differential form $\omega$. A non-trivial fact is that the discrete subgroup $K\subset\A$ is self orthogonal with respect to the pairing $\langle\cdot,\cdot\rangle$ (this corresponds to the Serre duality on $C$), ~$K^{\perp} = K$. It follows that for any element $g\in\A^*$, the annihilator of $gK$ is $g^{-1}K$, ~$(gK)^{\perp} = g^{-1}K$.

Let us choose an additive Haar measure $\mu$ on $\A$. The Fourier transform is defined as usual by the formula
$$
(\Fo f)(a)=\int_{\A}\langle a,b\rangle f(b)d\mu(b).
$$
Note that for any function $f\in\Dc(\A)$, the Fourier transform $\Fo (f)$ is well-defined and  belongs to $\Dc(\A)$.
Hence the Fourier transform $\Fo$ is well-defined on the dual space $\Dc'(\A)$ by the formula
$$
(\Fo\Delta)(f):=\Delta(\Fo^{-1}f)
$$
for all $\Delta\in \Dc'(\A)$ and $f\in \Dc(\A)$.
The Fourier transform  depends on the choice  of the differential form $\omega$. Namely, for $\Fo = \Fo_{\omega}$ we have
$$
(\Fo_{\omega'}f)(x) = (\Fo_{\omega}f)(ax), ~\mbox{if}~\omega' = a\omega.
$$
Furthermore, the Poisson summation formula and self-orthogonality of $K$ imply that
\begin{eqnarray}
\Fo\delta_K&=&\mu(\A/K)^{-1}\delta_K,\\\label{eq-pf3}
\quad \Fo\delta_{gK}&=&|g|^{-1}\mu(\A/K)^{-1}\delta_{g^{-1}K},\label{eq-pf4}
\end{eqnarray}
since $\mu(\A/gK)=|g|\mu(\A/K)$.

It is convenient to work with a self-dual measure $\mu$ on $\A$, which means by definition that for all functions $f,f'\in \Dc(\A)$, we have
$$
\Fo(\Fo f)(a)=f(-a),\quad\int_{\A}f(a)\overline{f'}(a)d\mu(a)=\int_{\A}\Fo f(a)\overline{\Fo f'}(a)d\mu(a),
$$
where the bar denotes the complex conjugation and the Fourier transform is taken with respect to $\mu$. Recall that a self-dual measure $\mu$ on $\A$ exists and is unique. Thus for a symmetric function $f$ and a symmetric distribution $\Delta$\footnote{Here, the word "symmetric" means that $(-1)^*f = f$ and $(-1)^*\Delta = \Delta$.}, we have $\Fo(\Fo f)=f$ and $\Fo(\Fo\Delta)=\Delta$. Also, it follows that $\mu(\A/K)=1$, and thus $F\delta_K=\delta_K$ and $F\delta_{gK}=|g|^{-1}\delta_{g^{-1}K}$. In addition, one can show that  $\mu(\OO)=q^{1-g(C)}$, where $g(C)$ is the genus of $C$.
 Locally, this means that  $ \mu (\OO_x) = q_x^{-1/2\nu_x(\omega)}$  and for a divisor $D$ on $C$, we have $\mu(\A(D))=q^{1-g(C)+\deg(D)}$. This implies that
\begin{eqnarray}
\Fo\delta_{K_x(D)}&=&q_x^{-1/2\nu_x(\omega) + \nu_x(D)}\delta_{K_x((\omega)-D)}\label{eq-Fourdiv-loc}\\
\Fo\delta_{\A(D)}&=&q^{1-g(C)+\deg(D)}\delta_{\A((\omega)-D)}\label{eq-Fourdiv-glob},
\end{eqnarray}
where $(\omega) = \sum_x \nu_x(\omega)\cdot x,~ D = \sum_x \nu_x(D)\cdot x, \quad \deg D =  \sum_x \nu_x(D)\cdot \deg(x)$.

\subsection{The Weil problem}

A simple but very important observation is that the topological group $\A^*$ acts
on the space $\A$ via multiplication. This defines
 representations of $\A^*$ on the complex vector spaces $\Dc(\A)$ and
$\Dc'(\A)$. Denote by an  upper index  the action of an element $g\in\A^*$ on a function or a distribution.
For all elements $g\in\A^*$, $a\in\A$, $f\in\Dc(\A)$, and $\Delta\in\Dc'(\A)$, we have
$$
f^g(a)=f(g^{-1}a),\quad \Delta^g(f)=\Delta(f^{g^{-1}}),
$$
$$
(-1)^*f(a)=f(-a),\quad (-1)^*\Delta(f)=\Delta((-1)^*f).
$$
In particular, for all elements $g,h\in\A^*$, we have
\begin{equation}\label{eq-actiondelta}
\delta_0^g=\delta_0,\quad \delta_{hK}^g=\delta_{g^{-1}hK},\quad \delta_{hK^*}^g=\delta_{g^{-1}hK^*}.
\end{equation}
Also, it is easy to show that for all $g\in\A^*$, $f\in\Dc(\A)$, and $\Delta\in\Dc'(\A)$, we have
\begin{equation}\label{eq-Fourrepr}
\Fo(f^g)=|g|(\Fo f)^{g^{-1}},\quad \Fo(\Delta^g)=|g|(\Fo\Delta)^{g^{-1}}.
\end{equation}

Since the group $\A^*$ is commutative, the
representation in $\Dc'(\A)$  decomposes  into one-dimensional
representations. The so-called Weil problem \cite{W1}\footnote{At the  same time, the problem how to describe homogenous generalized functions was considered and solved in \cite{GGPSh}[ch. 2].} is to find an eigen-distribution $\Delta_\chi\in\Dc'(\A)$ for each character $\chi$ as
above, that is  a distribution $\Delta_\chi$ such that for any element $g\in\A^*$, we have
$$
\Delta_\chi^g=\chi(g)\cdot\Delta_\chi.
$$
 Let  us  show that the distribution $\Delta_\chi$ is unique up to scalar. Indeed,  the space $\Dc(\A)$ is generated by linear combinations of the
 functions $\delta_{\A(D)} = \delta_{\OO}^g$ where $g \in \A^*$ such that $\A(D) = g^{-1}\A(0) = g^{-1}{\OO}$. Thus, the values of the distribution
 $\delta_{\chi}$ will be uniquely defined by its value on the function  $\delta_{\OO}$.
  Also, it follows that $\Delta_{\chi}$ is a symmetric distribution ($(-1)^*\Delta_\chi = \Delta_\chi$), provided it exists.

When ${\rm Re}(\chi)>1$, the eigen-distribution $\Delta_\chi$ can be explicitly expressed (up to scalar) as follows:
\begin{equation}\label{eq-int}
\Delta_\chi(f)=\int_{\A^*}\chi(g)f(g)d\mu^*(g),
\end{equation}
where $\mu^*$ is a multiplicative Haar measure on the locally compact group $\A^*$. It is useful to suppose that $\mu^*$ is normalized by the condition $\mu^*(\OO^*)=1$. Then $\mu^*$ is the  product of local Haar measures $\mu^*_x$ on locally compact groups $K_x^*$,  normalized by the conditions $\mu_x^*(\OO_x^*)=1$.

We can consider  a version of the Weil problem for the group  $\A^*$ acting on itself by multiplication. The solution  $\Delta^*_\chi$ will belong to the space $\Dc'( \A^*)$ and will be given by the integral $\eqref{eq-int}$ for a  function $f$ from the space $\Dc( \A^*)$, not the space $\Dc( \A)$. In that case, the integral will converge for all characters $\chi$. Thus, we get a new well-defined distribution
 $\Delta^*_\chi$.

 The original Weil problem can also be  reformulated in the following way. Let $i: \A^*\times \A \rightarrow \A$ be the multiplication map. It induces the direct image for function spaces
$$
i_*:\Dc( \A^*) \otimes\Dc( \A) \rightarrow \Dc( \A),\quad  i_*(f\otimes g) = f * g = \int \limits_{\A^*} f(x^{-1}y)g(y)\frac{d\mu(x)}{|x|}.
$$
This map is a an extension of the convolution map on the group  $\A^*$.
 And the original Weil problem is to find a distribution $\Delta_\chi \in \Dc'( \A)$  such that
 \begin{equation}\label{weil1}
 \Delta_\chi(f * g)) = \Delta^*_\chi(f)\cdot\Delta_\chi(g), \quad\mbox{in other words}\quad i^*\Delta_\chi =  \Delta^*_\chi\otimes \Delta_\chi.
  \end{equation}

\bigskip

Formula \eqref{eq-int} also clarifies  the relation between the Weil problem and the zeta-function of $C$. Namely, when ${\rm Re}(s)>1$, $s\in\C$, we have
$$
\Delta_{|\,\cdot\,|^s}(\delta_\OO)=\zeta_C(s),
$$
provided that $\mu^*$ is normalized as above. Indeed, by  \eqref{eq-int}, the left hand side equals
$$
\prod_{x\in C}\int_{K_x^*}|g_x|_x^s\delta_{\OO_x}(g_x)d\mu_x^*(g_x)=\prod_{x\in C}\left(\sum_{n\ge 0}q^{-ns\deg(x)}\right)=\prod_{x\in C}(1-q^{-s\deg(x)})^{-1},
$$
where $|\cdot|_x:=q^{-\deg(x)\nu_x(\cdot)}$\, and we have used that $\OO_x  =  \coprod_{n \geq 0}t_x^n\OO_x^*$ .

An analogous calculation shows that
\begin{equation}\label{eq-zetadiv}
\Delta_{|\,\cdot\,|^{s}}(\delta_{\A(D)})=q^{\deg(D)s}\zeta_C(s).
\end{equation}

Let us decompose the integral \eqref{eq-int} into an integral over the quotient $\A^*/K^*$
and sums over fibers of the projection $\A^*\to \A^*/K^*$:
\begin{equation}\label{eq-int2}
\Delta_\chi(f)=\int_{\A^*/K^*}\chi(g)\delta_{gK^*}(f)d\mu^*(g),
\end{equation}
where $f\in\Dc(\A)$, $g$ runs over all cosets in $\A^*/K^*$ and $\mu^*$ denotes also the induced measure on the quotient $\A^*/K^*$ (recall that the restriction of $\chi$ to $K^*$ is trivial). In what follows we will consider distribution valued integrals and thus  will remove  $f$ from the integral \eqref{eq-int2}.

The main difficulty with the Weil problem is that the formula \eqref{eq-int2} is not well-defined when ${\rm Re}(\chi)\le 1$, since the integral does not converge. However, the integral
\begin{equation}\label{eq-inthalf}
\int_{(\A^*)^{\Ge 1}/K^*}\chi(g)\delta_{gK^*}d\mu^*(g)
\end{equation}
over ``half'' of $\A^*/K^*$ is well-defined for any character $\chi$, where $(\A^*)^{\Ge 1}$ is the set of ideles $g\in \A^*$ with $|g|\ge 1$. Indeed, as above, the integral \eqref{eq-inthalf} converges when  ${\rm Re}(\chi)>1$. Besides, if $|g|\ge 1$, then $|\chi(g)|\le |\chi(g)||g|^N$ for any  $N>0$, and we have ${\rm Re}(\chi|\cdot|^N)>1$ for  sufficiently large $N$. Hence, one can bound  the integral \eqref{eq-inthalf} from above for any character $\chi$ by the convergent integral associated to  the character $\chi|\cdot|^s$ for a suitable~$s$.

\subsection{Regularization}

To overcome the convergence problem in the integral \eqref{eq-int2} one introduces a regularization of the integral.

Namely, let $\varphi$ be a function on $\A^*/K^*$ induced by a complex-valued function with finite support on $\Z$ via the degree homomorphism
$$
\deg:\A^*/K^*\to\Z.
$$
Thus, $\varphi$ has compact support and is locally constant on $\A^*/K^*$, in other words, $\varphi\in\Dc(\A^*/K^*)$. Define
$$
\Delta_\chi^\varphi:=\int_{\A^*/K^*}\varphi(g)\chi(g)\delta_{gK^*}d\mu^*(g).
$$
In can be easily seen that $\Delta^\varphi_\chi$ is well-defined for any character $\chi$ and that for any function $f\in\Dc(\A)$, the value of $\Delta^\varphi_\chi(f)$ is a Laurent polynomial in $q^{-s}$.

Consider the limit distribution $\lim\limits_{\varphi\to 1}\Delta_\chi^\varphi$, where $\varphi$ tends to the constant function with value $1$, which means that $\varphi(g)=1$ for all $g$ with $|\deg(g)|\le N$,  for increasing values of $N$. One can  prove directly that if ${\rm Re}(\chi)>1$, then this limit exists and is equal to $\Delta_{\chi}$. Now the idea is to show that the limit $\lim\limits_{\varphi\to 1}\Delta_\chi^\varphi$ exists for any $\chi$, is an eigen-distribution with respect to $\A^*$, and also satisfies a certain functional equation.

It is convenient to consider the decomposition $\varphi=\varphi_- +\varphi_+$\footnote{This corresponds to Riemann's decomposition of the domain of integration  of the integral defining the zeta-function $\zeta_C(s)$ as the  union of two pieces, $t < 1$ and $t > 1$, see above. More precisely, it corresponds  to  an approximation of  the  domain of integration by the two intervals  $N^{-1 }<  t < 1$ and $1 <  t < N$, with  $N \to \infty$ at the end of the computation.}, where for any $g\in\A^*/K^*$, we have $\varphi_-(g)=\varphi(g)$ if $\deg(g)<0$, $\varphi_-(g)=\varphi(g)/2$ if $\deg(g)=0$, and $\varphi_-(g)=0$ if $\deg(g)>0$. This  defines uniquely the second function $\varphi_+$.

We have
$$
\Delta^{\varphi}_\chi=\Delta^{\varphi_-}_\chi+\Delta^{\varphi_+}_\chi=
$$
$$
=\Delta^{\varphi_-}_\chi+
\int_{\A^*/K^*}\varphi_+(g)\chi(g)\delta_{gK}d\mu^*(g)
-\delta_0\cdot c(\varphi_+,\chi),
$$
where
$$
c(\varphi_+,\chi):=\int_{\A^*/K^*}\varphi_+(g)\chi(g)d\mu^*(g)\in\C.
$$
Suppose in addition that $\varphi$ is symmetric, that is, we have $\varphi(g)=\varphi(g^{-1})$. It follows that $\varphi_{\pm}(g)=\varphi_{\mp}(g^{-1})$. Using the Poisson summation formula (\ref{eq-pf4}) and  the decomposition \mbox{$g^{-1}K=g^{-1}K^*\cup\{0\}$}, we obtain
$$
\int_{\A^*/K^*}\varphi_+(g)\chi(g)\delta_{gK}d\mu^*(g)=
\int_{\A^*/K^*}\varphi_+(g)\chi(g)|g|^{-1}\Fo\delta_{g^{-1}K}d\mu^*(g)=
$$
$$
=\int_{\A^*/K^*}\varphi_+(g)\chi(g)|g|^{-1}\Fo\delta_{g^{-1}K^*}d\mu^*(g)+
(\Fo\delta_{0})\cdot c(\varphi_+,\breve\chi^{-1})=
$$
$$
=\int_{\A^*/K^*}\varphi_-(g)\chi(g)^{-1}|g|\Fo\delta_{gK^*}d\mu^*(g)+
(\Fo\delta_{0})\cdot c(\varphi_+,\breve\chi^{-1})=
\Fo\Delta^{\varphi_-}_{\breve\chi}+(\Fo\delta_{0})\cdot c(\varphi_+,\breve\chi^{-1}).
$$
In the second row, we replaced  $g$ by $g^{-1}$ and made use of the fact that $d\mu^*(g^{-1}) = d\mu^*(g)$.
Note that the Fourier transform is taken with respect to the self-duality of $\A$ provided  by a non-zero rational differential form $\omega$ and the self-dual additive measure $\mu$ on $\A$ (so that, $\mu(\A/K) = 1$).
Taking the sum, we get
\begin{equation}\label{eq-decompTI}
\Delta^{\varphi}_\chi=\Delta^{\varphi_-}_\chi+
\Fo\Delta^{\varphi_-}_{\breve\chi}+(\Fo\delta_{0})\cdot c(\varphi_+,\breve\chi^{-1})
-\delta_0\cdot c(\varphi_+,\chi).
\end{equation}
Since the integral \eqref{eq-inthalf} converges for any character $\chi$, the limits
$\Delta_\chi^-:=\lim\limits_{\varphi\to 1}\Delta_\chi^{\varphi_-}$ and
$\Fo\Delta_{\breve \chi}^-=\lim\limits_{\varphi\to 1}\Fo\Delta_{\breve\chi}^{\varphi_-}$ are well defined for any character
$\chi$.
In addition, for a given function $f  \in \Dc(\A)$,  the sequence $\Delta_\chi^{\varphi_-}(f)$  will stabilize when $\varphi\to 1$.

Furthermore, recall that the integral of a character on a compact commutative group equals zero if the character is non-trivial and equals the volume of the group if the character is trivial. This fact, together with a direct calculation, implies that $c(\varphi_+,\chi)=0$ if $\chi$ is not equal to $|\cdot|^s$ for some $s\in\C$, and that
$$
c(\varphi_+,|\cdot|^s)=\mu^*((\A^*)^{(1)}/K^*)\cdot\left(\sum_{n\in\Z}\varphi_+(g_n)q^{-ns}\right),
$$
where $g_n\in\A^*/K^*$ with $\deg(g_n)=n$. Therefore the limit
$c(\chi):=\lim\limits_{\varphi\to 1}c(\varphi_+,\chi)$ equals zero unless  $\chi = |\cdot|^s$ for some $s\in\C$, and
$$
c(|\cdot|^s)=\mu^*((\A^*)^{(1)}/K^*)\cdot\frac{1+q^{-s}}{1-q^{-s}}.
$$
It particular, $c(\chi)=-c(\chi^{-1})$.

By ~\eqref{eq-decompTI}, we obtain that the limit
$\Delta_\chi:=\lim\limits_{\varphi\to 1}\Delta^\varphi_\chi$ is well defined and equals
$$
\Delta^-_\chi+\Fo\Delta^-_{\breve \chi}+(\Fo\delta_0)\cdot c(\breve\chi^{-1})-\delta_0\cdot c(\chi).
$$
 Making use of formulas~\eqref{eq-actiondelta} and \eqref{eq-Fourrepr}, we see that $\Delta_\chi$ is indeed an eigen-distribution with the character $\chi$. In addition, we obtain  the functional equation for distributions
$$
\Delta_\chi=\Fo\Delta_{\breve\chi}
\footnote{This equation shows that the RHS does not depend on the choice of the differential form $\omega$, which enters in the definition of   the Fourier transform. This fact   follows easily  from the behavior of the distribution under translation and from the triviality of the character $\chi$ on the subgroup $K^*$ . }.
$$
Note that up to scalar this equation also follows from formulas~\eqref{eq-actiondelta} and \eqref{eq-Fourrepr}, and the fact that $\Delta_\chi$ is an
eigen-distribution. Thus we have shown that the Weil problem has a solution for each character $\chi$.

\bigskip

\subsection{Applications}
This has the following application to the zeta-function of $C$.
First, since \mbox{$\Delta_{|\,\cdot\,|^s}(\delta_\OO)=\zeta_C(s)$} for ${\rm Re}(s)>1$, we deduce that $\zeta_C(s)$ has a meromorphic continuation on the whole $s$-plane. In addition, the formulas~\eqref{eq-Fourdiv-glob} and \eqref{eq-zetadiv} imply that the functional equation
$$
\zeta_C(s)=q^{g(C)-1}q^{(2-2g(C))(s)}\zeta_C(1-s)
$$
is satisfied, since  $\deg((\omega))=2g(C)-2$.
Finally, it is easy to see that
$$
\mu^*((\A^*)^{(1)}/K^*)=\frac{|\Pic^0(C)(\Fb_q)|}{q-1},
$$
provided that $\mu^*(\OO^*)=1$, where $\Pic^0(C)$ is the group of isomorphism classes of degree zero line bundles on $C$ and $\Pic^0(C)(\Fb_q)$ its  group of rational points. Indeed, one has the exact sequence
\begin{equation}\label{eq-exseq}
1  \rightarrow \Fb_q^*  \rightarrow \OO^*  \rightarrow (\A^*)^{(1)}/K^* \rightarrow \Pic^0(C)(\Fb_q)  \rightarrow 1,
\end{equation}
the group  $\OO^*$ is compact (a trivial assertion) and
 the group $\Pic^0(C)(\Fb_q)$ is finite (a much less trivial one). Therefore, the only  singularities of $\zeta_C(s)$ are two first order poles at $s=0$ and $s=1$ with  residues
$$
\res_{s = 0} \zeta_C(s) =  -\frac{|\Pic^0(C)(\Fb_q)|}{q-1}\ln^{-1}q
$$
and
$$
\res_{s = 1} \zeta_C(s) = q^{1-g(C)}\frac{|\Pic^0(C)|(\Fb_q)}{q-1}\ln^{-1}q.
$$
In addition, we can easily deduce the following presentation for $\zeta_C(s)$
$$
\zeta_C(s) =   \frac{P(q^{-s})}{(1 - q^{-s})(1 - q^{1-s})}
$$
where $P(q^{-s})$ is a polynomial in the variable $t =  q^{-s}$ of degree $2g(C)$ (see details in \cite{W2}[ch. VII, \S 6, Theorem 4]).
This polynomial satisfies the functional equation $P(q^{-1}t^{-1}) = q^{-g(C)}t^{-2g(C)}P(t)$. Also, we can remove the singularities
and get for $\zeta_C(s)$ the presentation
$$
 -\frac{|\Pic^0(C)(\Fb_q)|}{q-1}
\frac{1}{1 - q^{-s}}  +  q^{1-g(C)}\frac{|\Pic^0(C)|(\Fb_q)}{q-1}\frac{1}{1 - q^{1-s}}  + \int \limits_{\A^*_{-}} [|g|^sf(g)   +  |g|^{1-s}\Fo^{-1}f(g)] d\mu^*(g).
$$
where the integral over the "domain" $\A^*_{-}$ means the sum of integrals over  $\{ g \in \A^*:  |g | > 1\}$ and $1/2$ of the integral
 over $\{ g \in \A^*:  |g | = 1\}$. The functions under the integral will be $f(g) = \delta_{\OO}$ and $\Fo^{-1}f(g) = \Fo f(g) = q^{1-g(C)}\delta_{\A((\omega))}$.

 To see the analogy with  the case of Riemann's zeta function (\ref{eq-rz}) we note that $\ln q/(1 - q^{-s}) \sim  1/s$ and  $\ln q/(1 - q^{1-s}) \sim  1/(s-1)$. We will refer to  the presentation which we obtain as the principal parts decomposition of the meromorphic function $\zeta_C(s)$.

\bigskip

In our discussion, we implicitly supposed that the group $(\A^*)^{(1)}/K^*$ is compact. This  fact is a non-trivial one. It includes both finiteness
of the class number (for number fields and  $\Pic^0(C)(\Fb_q)$ in our case) and the Dirichlet theorem on units (see \cite{A, W2}).

It is remarkable and very important that these finiteness results are consequences of the analytical  computations we have already carried out.
Namely, we can apply  both parts of  (\ref{eq-decompTI}) to  a given function $f$, set $\chi = |\cdot|^s$  and then the  decomposition
(\ref{eq-decompTI}) can be rewritten as
$$
(\mbox{polynomial in}~q^{-s}) =  (\mbox{polynomial in}~q^{-s}) - \delta_{(0)}(f)c(\varphi_+, |\cdot|^s) + \delta_{(0)}(Ff)c(\varphi_+, |\cdot|^{1-s}).
$$
Choosing the function $f$ in an appropriate way, we get
$$
(\mbox{polynomial of}~q^{-s}) =  (\mbox{polynomial of}~q^{-s}) + (\mbox{nonzero polynomial of}~q^{-s} ) \int \limits_{(\A^*)^{(1)}/K^*}1.
$$
Thus, the  volume of  $(\A^*)^{(1)}/K^*$ is finite and the group itself is compact (see above, the exact sequence (\ref{eq-exseq}) for this group).
This observation was made by Iwasawa \cite{I2}.

{\bf Scholium 1}. An important  remark is that the entire analytical structure (``singularities" and  ``principal parts") and the symmetry properties of zeta- and $L$-functions can be visualized without  analytic continuation.
They are  present in their entirety even for the regularized zeta-function {\bf  before}  going to the final limit\footnote{That is, $\varphi \to 1$.}.

\section{Discrete versions and holomorphic tori}
 \subsection{Analysis on discrete groups}
 Let $C/\mathbb F_q$~ be again an algebraic curve defined over a finite field $\mathbb F_q$. We have on $C$ the following (discrete) groups of divisors:
$$\Gamma = \Div(C) \supset \Gamma_0 = {\Div}^0 (C) \supset \Gamma_l = {\Div}_l (C),
$$
where  $\Div^0 (C)$ is the group of divisors of degree 0 and $\Div_l (C)$ is the group of  divisors linearly equivalent to 0.
As well known, $ \Gamma_0/ \Gamma_l$ is  a finite group $\Pic^0(C)(\Fb_q)$.

Before the Tate-Iwasawa approach had been developed, another method was known for the case of algebraic curves. It was based on a direct application of the Riemann-Roch theorem (see, for example, \cite{D}). Denote by $\Gamma_+ = \{ D \in \Gamma : D \ge 0\}$ the subset of effective divisors on $C$. One can write the zeta-function as the Dirichlet series  converging  for ${\rm Re}(s)  > 1$
 $$
\zeta_C(s) = \sum_{D \in \Gamma_+} \Nm(D)^{-s}, \quad  \mbox{where}~\Nm(D) = q^{\deg D}.
$$
Assuming for simplicity that there is a point $P \in C$ such that $k(P) = \Fb_q$ and choosing  a finite set of representatives $R$ of cosets of
$ \Gamma_l$ in  $\Gamma_0$,  we have
the following decompositions:
$$
\Gamma = \coprod_{R, n \in \Z}(R + nP)\Gamma_l, \quad \Gamma_+ =  \coprod_{R, n \in \Z}(R + nP)\Gamma_l \cap \Gamma_+.
$$
Therefore, the Dirichlet series can be transformed as\footnote{Here we use the standard notations: $\sim$ for the linear equivalence of divisors and $l(D)$
for  dimension of  the space $\mbox{H}^0(C, {\cal O}_C(D))$. }
$$
\zeta_C(s) =   \sum_{R, n \in \Z} \#\{D \in \Gamma : D \sim  R + nP, ~ D \ge 0 \}q^{-ns} =
$$
$$
 \sum_{R, n \ge 0}\frac{q^{l(R + nP)} - 1}{q - 1} q^{-ns} = - \frac{1}{q - 1}\sum_{R, n \ge 0} q^{-ns} + \frac{1}{q - 1} \sum_{R, n \ge 0}q^{l(R + nP)} q^{-ns} = ?
 $$
Now, we come to the main point, which  is to use the duality and Riemann-Roch together: $ l(R + nP) = l({\rm K} - R - nP) + n + 1 - g(C)$. Here, ${\rm K} = (\omega),~\omega \in \Omega^1_K$ is the canonical
class and   we choose  $\omega \in \Omega^1_K$ such that ${\rm K} = (2g -2)P + R_0$ with $g = g(C),  R_0 \in \{R\}$. Note that $\{R_0 - R\}  =  \{S \}$ where $S$ is another finite set of representatives  of cosets of
$ \Gamma_l$ in  $\Gamma_0$. We continue the transformation of the zeta-series:
$$
=  - \frac{1}{q - 1}\sum_{R, n \ge 0} q^{-ns} + \frac{q^{1 - g}}{q - 1} \sum_{R, n \ge 0}q^{l(R_0 - R +(2g -2 - n)P)} q^{-n(s-1)} =
$$
$$
- \frac{1}{q - 1}\sum_{R, n \ge 0} q^{-ns} + \frac{q^{1 - g}}{q - 1}\sum_{S, n \ge 0} q^{-n(s-1)} +  \frac{q^{1 - g}}{q - 1}\sum_{S, 2g-2 \ge n \ge 0}(q^{l( R +(2g -2 - n)P)}-1) q^{-n(s-1)}.
$$
At the end, we get sum of three terms. The first one has poles at the points of arithmetical progressions $2\pi im/\ln q,~m \in \Z$, the second has poles at the points of arithmetical progressions $1 + 2\pi im/\ln q,~m \in \Z$ and the last one is  a polynomial in $q^{-s}$.
Thus, we made an analytic continuation of  $\zeta_C(s)$ to the whole $s$-plane. It is also not difficult to obtain  the functional equation \cite{D}[\S 23].
It is worthwhile to compare this computations with the ones made in the previous section.
Elaborating this approach in the spirit of the Tate-Iwasawa  method, one can also incorporate a function $f = f(D)$
into the Dirichlet series and study more general sums
$$
\zeta_C(s, f) = \sum_{D \in \Gamma}f(D) \Nm(D)^{-s}.
$$
We  leave this to the reader but will now  introduce and carefully study the necessary tools: function spaces, the Fourier transform and the Poisson
formula.

\bigskip

We  filter the group $\Gamma$ by the subgroups  $\Gamma_S : =\Div_S(C) =$
 \{divisors with support in $S$\}:
$$\Gamma = \varinjlim_S\Gamma_S = \varinjlim_S\bigoplus_{x\in S}  \Gamma_{(x)} = \bigoplus_{x\in C} \Gamma_{(x)},$$
 where the $S$ are  finite subsets of  $ C$ and the $\Gamma_{(x)}= \Z$ is the valuation group at the point $x$.

On $\Z$ there are two important spaces: the space ${\cal D}(\Z)$ of finitely supported functions and the space ${\cal D}_+(\Z)$ of functions $f$ such that $f(n) = 0~ \mbox{for}~n \ll 0  ~\mbox{and}~f(n) = \mbox{const}~ \mbox{for}~n \gg 0$.

 The space ${\cal D}_+(\Z)$ is  a linear span of the functions
$$
\delta_{(\ge m)}(n) = \left\{
	\begin{array}{ll}
		1,& n\geqslant m\\
		0,& n <  m .	
	\end{array}
\right.
$$
The function space which we've defined can be directly deduced from the Bruhat spaces which we introduced above.
The group $\Z = \Gamma_{(x)}$ appears as $K_x^*/ \hat{\OO}_x^*$ where  $K_x$ is the local field at the point $x \in C$. Then
$K_x = \coprod_{n \in \Z} t_x^n \hat{\OO}_x^* \cup \{0\}$.
Thus, we have the isomorphism
\begin{equation}\label{isom}
{\cal D}(K_x)^{\hat{\OO}_x^*} \stackrel{\sim}\rightarrow {\cal D_+}(\Z), \quad \{f \mapsto f(t_x^n)\},
\end{equation}
from  the  subspace of functions invariant under multiplication by $\hat{\OO}_x^*$ onto the space of functions on the group $\Z$.
In the same way, we have an isomorphism
$$
{\cal D}(K_x^*)^{\OO_x^*} \stackrel{\sim}\rightarrow {\cal D(\Z)},
$$
where ${\cal D}(\Z)$ is the space of functions on $\Z$ with finite support.

Just as in the case of Bruhat functions, we have the convolution map
$$
 * : {\cal D}(\Z) \otimes {\cal D}_+(\Z) \rightarrow  {\cal D}_+(\Z),
$$
which sends $f \otimes g$ to the function $f * g$ defined by
$$
f*g(n) = \sum_{m \in \Z}f(n - m)g(m),\quad f \in {\cal D}(\Z), ~g \in {\cal D}_+(\Z)
$$
and the space  ${\cal D}_+(\Z)$ is a free ${\cal D}(\Z)$-module with generator $
\delta_{(\ge 0)}$.

\bigskip

The Fourier transform   $\Fo = \Fo_x: {\cal D}(K_x) \rightarrow {\cal D}(K_x)$  commutes with
the action of the group $\hat{\OO}_x^*$ and consequently  induces an analogue of the  Fourier transform
on the group $\Z = \Gamma_{(x)}$ at the point $x$:
$$
\Fo_x: {\cal D_+}(\Z) \rightarrow {\cal D_+}(\Z).
$$
Recall that we use a self-dual measure on the local field $K_x$ to define the Fourier transform $\Fo_x$. Let $k_x : = \nu_x(\omega)$ for $\omega \in \Omega^1_K$.Then, one can easily compute this map on the standard functions
$$
\Fo_x (\delta_{(\ge m)})   =  q_x^{-1/2k_x - m}\delta_{(\ge -k_x - m)}.
$$
Here we use  that if $\nu_x(D) =  n_x$ then $\delta_{K_x(D)}$ is mapped to $ \delta_{\ge -n_x}$ under the morphism (\ref{isom}) (see section 2.2).

\bigskip

Let us now globalize these purely local considerations. We need spaces  ${\cal D}(\Gamma)$ and  ${\cal D}_+(\Gamma)$ for
the group of all divisors. First, we change the notations for the standard functions:
$$
f_{x,D}(n_x) = \left\{
	\begin{array}{ll}
		1,& n_x\geqslant \nu_x(D)\\
		0,& n_x < \nu_x(D),	
	\end{array}
\right.
$$
where $n_x \in \Gamma_{(x)}$ and  $D$ runs through the divisors on the curve $C$. We will sometimes denote this function by the old notation $\delta_{(\ge n)}$ where $~n = \nu_x(D)$.

Taking finite products
$$
f_D((n_x)_{x \in S}) = \otimes_{x\in S} f_{x, D}(n_x),
$$
we get  functions on $\Gamma_S$.
We set
$$
{\cal D}(\Gamma):  =  \{ \mbox{linear span of the functions}~f_D = \otimes_{x \in \Supp D}f_ {x, D}~\otimes_{x \not\in \Supp D}\delta_{(0)} \},
$$
where $\delta_{(0)}(n)  =  1$ if $n = 0$ and $= 0$ if $n \neq 0$. This  space consists of all finitely supported functions on  $\Gamma$.
  Also, we have the space
$$
{\cal D}_+(\Gamma):  = \{ \mbox{linear span of the functions}~f_D = \otimes_{x \in \Supp D}f_ {x, D}~\otimes_{x \not\in \Supp D}\delta_{(\ge 0)} \}
$$
which is equal  to  the space of invariants ${\cal D}(\A^*)^{\OO^*}$
where $\OO^* = \prod_{x \in C}\OO^*_x$. The functions $f_D$  correspond to the functions $\delta_{\A(-D)} \in {\cal D}(\A)^{\OO^*}$.

\bigskip

The Poisson formula $\delta_{gK}(f) = |g|^{-1}\delta_{g^{-1}K}(\Fo f)$ for the space ${\cal D }(\A)$ can be formulated in the new situation as follows
\footnote{Since $- 1 \in {\OO}^*$ we have $\Fo = \Fo^{-1}$ on the spaces ${\cal D }(\Gamma)$ and ${\cal D }_+(\Gamma)$.}.
For simplicity, we first suppose that $g = 1$. Rewrite  the Poisson formula as
$$
\delta_{(0)}(f) + \delta_{K^*}(f) = \delta_{(0)}(Ff) + \delta_{K^*}(Ff)
$$
and assume that the function  $f$ belongs to the space   ${\cal D}(\A)^{\OO^*}$. The group $\Gamma$ contains  the subgroup  $\Gamma_l$,  which is   the group of divisors of functions $g \in K^*$. Thus, $\Gamma_l = K^*/\Fb_q^*$. This  means that the delta function  $\delta_{K^*}(f)$ is equal to  $(q-1)$ times the delta-function $\delta_{\Gamma_l}(f)$,  applied to $f$ considered as a function on  the space ${\cal D}_+(\Gamma)$.
What should we  do with the functionals  $\delta_{(0)}(f)$ ? The zero subgroup in $\A$ does not correspond to any subset in  $\Gamma$ and the functional can
be defined  in the following way:
$$
\delta_{(0)}(f) = \prod_x  \delta_{(0)}(f_x),~\mbox{if}~f = \otimes_{x \in C} f_x,
$$
with
$$
\delta_{(0)}(f_x) = \lim_{n \to \infty}f_x(n).
$$
The new Poisson formula for the functions $f \in {\cal D}_+(\Gamma)$  reads
\begin{equation}\label{eq-pf1}
\delta_{(0)}(f)  + (q-1)\delta_{\Gamma_l}(f) = \delta_{(0)}(\Fo f)  + (q-1)\delta_{\Gamma_l}(\Fo f).
\end{equation}
More generally, when $g \neq 1,~g \in \A^*$, we get
\begin{equation}\label{eq-pf5}
\delta_{(0)}(f)  + (q-1)\delta_{\gamma+\Gamma_l}(f) = \delta_{(0)}(\Fo f)  + (q-1)\delta_{-\gamma+\Gamma_l}(\Fo f)
\end{equation}
with $\gamma \in \Gamma$ (which is the image of $g$ under the natural map $\A^* \rightarrow \Gamma$).

\bigskip

\subsection{Extended group of divisors}
The presence of the multipliers $q - 1$ can be avoided if we consider an extended adelic group.  Let
$\mathcal{O}_x/\mathfrak{m}_x = k(x)$ and let us set
$$
\tilde\Gamma_x: = K^*_x/(1+\mathfrak{m}_x).
$$
We then have  a canonical exact sequence:
$$
1\longrightarrow k(x)^* \longrightarrow \tilde\Gamma_x
\stackrel{\pi}{\longrightarrow} \mathbb{Z}\longrightarrow 0.
$$
For a subset $M \subset \Z$ we define $\tilde\Gamma_x(M):=\pi^{-1}(M)$.
Then we introduce the following space of functions on this new group
$$
{\cal D}_+(\tilde\Gamma_x)=
\left\{
\begin{array}{llllll} f\in \mathcal{F}(\tilde\Gamma_x):\space &
f(\gamma) & = & \mbox{const} & \mbox{\space on \space}\tilde \Gamma_x\mbox{\space}
(\geqslant n)\mbox{\space} & \mbox{\space for\space} n\geqslant n_1\\
 & f(\gamma) & = & 0 & \mbox{\space on \space} \tilde\Gamma_x \mbox{\space}
(n) \mbox{\space} & \mbox{\space for\space}  n \leqslant n_0
\end{array}
\right\}
$$
 We can then introduce the global group
$\tilde\Gamma_C=\prod \nolimits^{'}_{x\in C}\tilde\Gamma_x$, which  fits into
the exact sequence
$$\begin{array}{ccccccccc} 1 & \longrightarrow  & \prod
\limits_{x\in C}k(x)^*  & \longrightarrow  & \tilde\Gamma_C &
\longrightarrow  & \bigoplus\limits_{x\in C} \mathbb{Z} &
\longrightarrow  & 0\\
&&&&&&\parallel\\
&&&&&&\mathrm{Div}(C) \end{array}$$
As above, fixing a differential form  $\omega\in \Omega^1_{\mathbb{F}_q(C)}~\omega\ne0$,
$(\omega)=\sum k_x\cdot x$,  we define the local Fourier transforms
$\Fo: \mbox{\space} {\cal D}_+(\tilde\Gamma_x)\overset \sim \to {\cal D}_+(\tilde\Gamma_x)$
and its global counterpart.
There is  a map
 $K^*=\mathbb{F}_q(C)^*\rightarrow \tilde\Gamma_C$ and if we set
$\tilde\Gamma_l$ being the image of $K^*$  under this map, we have a new
Poisson formula (without additional multipliers $(q - 1)$):
$$
\boxed{\delta_{(0)}(f)+\delta_{\gamma+\tilde\Gamma_l}(f)= q^{- \deg(\gamma)}\delta_{(0)}(\Fo f)+q^{- \deg(\gamma)}\delta_{-\gamma+\tilde\Gamma_l}(\Fo f)}
$$
with $\gamma\in\tilde\Gamma/\tilde\Gamma_l$.
Here, we can define the functionals $\delta_{(0)}(f)$ as we did  for the group $\Gamma$.

Indeed, inside the full adelic group we have subgroups $K^*\!\subset \!\mathbb{A}^*$ and  $\prod \limits_x
(1+\mathfrak{m}_x)=:\mathcal{O}^*_1 \!\subset \!\mathbb{A}^*$.
Then  $\tilde\Gamma_C=\mathbb{A}^*/\mathcal{O}^*_1$ and $\mbox{\space}
K^*\cap\mathcal{O}^*_1 = 1$. This  gives Poisson on $\mathbb{A}_K$ for functions $f\in
{\cal D}(\mathbb{A})^{\mathcal{O}^*_1}={\cal D}_+(\tilde\Gamma_C)$.

There is a measure on $\tilde\Gamma_C$,  by counting on the discrete group $\Gamma_C = \Div(C)$
and  integration over the compact group $\prod\limits_x k(x)^*$.  We can then make the entire computation of
the zeta-function $\zeta_C(s)$ using only analysis  on the group  $\tilde\Gamma_C$
\footnote{ On the curve $C/\mathbb{F}_q$, there is a Heisenberg group
$1\longrightarrow \bigoplus \limits_{x\in C}k(x)^* \longrightarrow
\hat{\Gamma}_C \longrightarrow {\tilde\Gamma}_C  \longrightarrow 1$
or
$$
\begin{pmatrix} 1 & \bigoplus \nolimits_x \mathbb{Z} &
\bigoplus \nolimits_x
k(x)^{*}\\
0 & 1 & \prod \nolimits_x k(x)^{*}\\
0 & 0 & 1 \end{pmatrix}
$$
if we put
$(n_x)\times(b_x)\longmapsto (b_x^{n_x})$ for  $(n_x)\in\bigoplus\nolimits_x
\mathbb{Z}$, $(b_x)\in\prod\nolimits_x k(x)$, $(b_x^{n_x})\in\bigoplus\nolimits_x
k(x)$. What is its arithmetic meaning ? Are there some $L$-functions related to this group and, in particular,
to characters of its center?}.

\bigskip

\subsection{Holomorphic tori}
The next step in  our construction will be  to consider the dual groups of the groups such as $\Z = \Gamma_{(x)}$ or $\Gamma_C$. The usual Pontryagin definition is not adequate for our arithmetical goals. We wish  to consider the discrete groups as commutative group schemes over $\C$. Inside the category of commutative group schemes we have two subcategories, a category ${\cal T}$ of commutative algebraic groups of toroidal type  (such that the connected components are  tori isomorphic to several copies of ${\mathbb G}_m$) and a  category ${\cal A}$ of finitely generated abelian groups (considered as the discrete group schemes). We have two functors between them:
$$
\Theta : {\cal T}\rightarrow{\cal A}\quad \Phi(T) = \Hom(T, {\mathbb G}_m),
$$
$$
\Psi : {\cal A} \rightarrow{\cal T}\quad \Psi(A) = \Hom(A, {\mathbb G}_m)
$$
with the obvious definition of the corresponding maps on arrows.
If $G, G' \in \mbox{Ob} {\cal A}$ and $f: G \rightarrow G'$ is an arrow of ${\cal A}$, we denote the dual map
$$
\check f: \Psi(G') \rightarrow \Psi(G).
$$
These functors satisfy by all the standard properties for the duality maps (see [{\bf SGA 3}]).

The dual groups of our discrete groups will be complex  tori:
$$
{\mathbb T}_{(x)}: = \Hom(\Gamma_{(x)}, \C^*)
$$
$$
{\mathbb T}_{\Gamma}(\mbox{or, simply}~{\mathbb T})  = \Hom(\Gamma, \mathbb C^*):  = \varprojlim_S\Hom(\Gamma_S,\mathbb C^*),
$$
$$
\mathbb T_{\Pic} := \Hom(\Gamma/\Gamma_l,\mathbb C^*) \supset  \mathbb T_0 := \Hom(\Gamma /\Gamma_0, \mathbb C^*) = \mathbb C^*,
$$
where the inclusion is induced by the degree map: $\Gamma \stackrel{\deg}{\rightarrow} \Gamma / \Gamma_0 = \mathbb Z$\footnote{The complex structure on the set characters (with values in $\C^*$) was considered by A. Weil in \cite{W2}.}.

We  have a short exact sequence
$$1 \rightarrow \mathbb T_0 \rightarrow \mathbb T_{\Pic} \rightarrow  \Hom(\Div^0(C)/\Div_l(C),\mathbb C^*) \rightarrow 0$$
and therefore,
$\mathbb T_{\Pic} \cong \mathbb C^*\times \{\mbox{finite group}\}$. In the sequel, . we will denote the finite group $\Div^0(C)/\Div_l(C)$ by $\Phi$ and its dual group as $\check\Phi$.

In order to consider the torus $\mathbb T_{\Gamma}$ as a manifold one only needs to know the ring of regular functions on $\mathbb T_{\Gamma}$.
For this  we introduce the projective system of  tori  $\mathbb T_S:=\Hom(\Div_S(C),\mathbb C^*) = \prod_{x \in S}\Hom (\Gamma_{(x)}, \C^*) = :
\prod_{x \in S}{\mathbb T}_{(x)}$ with a natural projection
$$\mathbb T_S \leftarrow \mathbb T_{S'}, \quad\mbox{if}~S\subset S',$$
induced by embedding $\Div_S(C)\subset \Div_{S'}(C)$.

Let us now consider the divisor $D_S$ with normal crossings  on  ${\mathbb T}_S$ that consists of the points in the product ${\mathbb T}_S$  for which  at least one component is the identity point in some ${\mathbb T}_{(x)}$. Let $\C_+[{\mathbb T}_S]$ be the space of  rational functions on ${\mathbb T}_S$ that are regular outside  $D_S$ and may have poles of  first order on $D_S$. In particular,  $\C_+[{\mathbb T}_{(x)}] = \C[{\mathbb T}_{(x)}](1 - z)^{-1}$.

Then $\mathbb T: = \varprojlim_S \mathbb T_S$ and we set
$$
\mathbb C[\mathbb T] := \varinjlim_S \mathbb C [\mathbb T_S]  = \bigotimes\nolimits_{x\in C}^{\prime} \mathbb C[\mathbb T_{(x)}]: = \{(\otimes_{x\in C}f_x)\colon f_x \in \mathbb C[\mathbb T_{(x)}]~  \& ~f_x = 1~\mbox{for almost all}~ x\}
$$
$$
\mathbb C_+[\mathbb T] := \varinjlim_S \mathbb C_+ [\mathbb T_S] =  \bigotimes\nolimits_{x\in C}^{\prime} \mathbb C_+[\mathbb T_{(x)}]: = \{ (\otimes_{x\in C} f_x)\colon f_x \in \mathbb C_+[\mathbb T_{(x)}]~ \& ~f_x = (1 - z_x)^{-1}~\mbox{for almost all}~ x\}.
$$
Here  $z_x$ is a coordinate on $\mathbb T_{(x)}$, namely, for $\chi \in \mathbb T_{(x)}$ we have $z_x = \chi(1)$.

We can easily see that
\begin{enumerate}
	\item $\mathbb C_+[\mathbb T]$ is a free module over the ring  $\mathbb C[\mathbb T]$ with generator $(f_x)$ where
$ f_x = (1 - z_x)^{-1}~\mbox{for all}~ x$.
	\item  for all (finite) $S$, $\mathbb C_+[\mathbb T_S]\supset\mathbb C[{\mathbb T}_S]$ but nevertheless $\mathbb C_+[\mathbb T]\not\supset
\mathbb C[\mathbb T]$.
\end{enumerate}

Let us now assume that we have an algebraic group  $\mathbb T$ which is not necessarily connected but whose  connected component  $\mathbb T_e$ of the identity  $e$ is  a torus,  and for which
$\pi_0(\mathbb T) = \mathbb T/\mathbb T_e$ is a discrete group. Then we know that
$$
\mathbb C[\mathbb T] = \bigoplus_{h \in \pi_0(\mathbb T)}\mathbb C[h\mathbb T],
$$
and we set
$$
\mathbb C_+[\mathbb T] = \mathbb C_+[\mathbb T_e]\bigoplus_{h \in \pi_0(\mathbb T),~h \ne e}\mathbb C[h\mathbb T].
$$

 We would now like to describe the  relation between the functions on the discrete groups and on  their duals. Certainly, this must be some kind of  Fourier transform but to distinguish it from the Fourier transform on the original discrete groups we will
call the corresponding transformation  the Mellin transform\footnote{Certainly, this usage corresponds to the well known cases where the classical Mellin map has appeared.}.
If $f  \in {\cal F}(\Z)$ belongs to a function space on $\Z$ then the Mellin transform ${\rm M}f$ will be
\begin{equation} \label{l-tr}
{\rm M}f(z) = \sum_{n \in \Z} f(n)z^n,
\end{equation}
where $z \in \Hom (\Z, \C^*)$ and $z^n$ is $z(n)$.
The Mellin transform respectively  maps the spaces ${\cal D}(\Z)$ and ${\cal D}_+(\Z)$ isomorphically onto the spaces $\mathbb C[{\mathbb T}_{(x)}]$ and $\mathbb C_+[{\mathbb T}_{(x)}]$. These isomorphisms allow us to define the  Fourier transform, say  $\Fo_{\mbox{new}}$, on the spaces $\mathbb C_+[{\mathbb T}_{(x)}]$ such that ${\rm M}\circ\Fo_{\mbox{old}}=\Fo_{\mbox{new}}\circ{\rm M}$\footnote{Here, $\Fo_{\mbox{old}}$ is the Fourier transform on the space  ${\cal D}_+(\Z)$. In the sequel, we will use the same notation  $\Fo$ for the  Fourier transforms in all spaces.}.
The Fourier transforms of the standard functions on $\mathbb T_{(x)}$
can be computed in an explicit way. Starting from the equality $$
\Fo_x(f_{x,D}) = q_x^{-\frac12 k_x + n_x} f_{x, {\rm K}-D}
$$
on the group $\Gamma_x$,
where $D = \sum_x n_x\cdot x$,~$n_x = \nu_x(D)$ and ${\rm K} = \sum_x k_x\cdot x$~is the canonical divisor,
we get in the  space  $\C_+[\mathbb T_{(x)}] = \mathbb C[z_x, z_x^{-1}](1-z_x)^{-1}$ for any $m \in \Z$
$$
\Fo_x(z_x^{m}(1-z_x)^{-1}) = q_x^{-\frac12 k_x - m}z_x^{-(k_x + m)} (1-z_x)^{-1}
$$
since
$${\rm M}f_{x,D} = z_x^{\nu_x(D)} (1-z_x)^{-1} = z_x^{n_x} (1-z_x)^{-1}.$$
From the definition  for the standard functions we can find
 the Fourier transform in the general form,  for an  arbitrary function $f$ from  $\C_+[\mathbb T_{(x)}]$
$$\Fo_x(f(z_x)) = q_x^{-\frac12 k_x} z_x^{-k_x} (1-q_x^{-1}z_x^{-1}) (1-z_x)^{-1} f(q_x^{-1}z_x^{-1}). $$

 We now move to the global construction.
Applying  the Mellin transform to the standard functions, we check that
$$ \C_+[\mathbb T_S] = \mbox{linear span of functions}~{\cal L}f_D, D\in\Div_S(C).$$
We also see that  for almost all $x \in C$,  $\Fo_x((1 - z_x)^{-1}) = (1 - z_x)^{-1}$ in the space $\mathbb C[{\mathbb T}_x]$.
These facts allows us  to combine   the local Fourier transforms $\Fo_x$ for all $x$ into a global one $\Fo$ on the torus ${\mathbb T}_{\Gamma}$
$$
\Fo\colon \C_+[{\mathbb T}_{\Gamma}] \stackrel{\sim}{\rightarrow} \C_+[{\mathbb T}_{\Gamma}],
$$
 such that $\Fo(f_x) = (\Fo_x f_x)$ and  $\Fo\circ \Fo = \Id$.

\bigskip

The exact sequence
$$0 \rightarrow  \Gamma_0 \rightarrow  \Gamma  \stackrel{\pi}\rightarrow \Gamma/\Gamma_0 \rightarrow 0$$
determines an embedding of  tori:
$$
\check\pi: \mathbb T_0 \rightarrow  \mathbb T.
$$
If $z_x\in\mathbb T_x = \Hom(\mathbb Z_{(x)},\mathbb C^*)$ and $\Div(C) = \bigoplus_{x\in C}\mathbb Z_{(x)}$, then the embedding $\check\pi\colon\mathbb T_0\to\mathbb T$ can be written as  $z\mapsto z_x = z^{\deg(x)}$ where $\#k(x) =q_x=q^{\deg(x)}$.

The embedding  induces a  map of  $\C_+[{\mathbb T}_{\Gamma}]$ into some function space
\begin{equation}\label{eq-emb}
\C_+[\mathbb T]\xrightarrow{\check\pi^*}{\cal F}(\mathbb T_0),
\end{equation}
where we still have  to identify the target space !
By direct substitution, one gets
$$
\check\pi^*({\cal L}f_D) = \prod_{x\in C}\left(z^{\deg(x)}\right)^{\nu_x(D)} \left(1-z^{\deg(x)}\right)^{-1} =
$$
$$
= z^{\deg D}\prod_x \left(1-z^{\deg(x)}\right)^{-1} = z^{\deg(D)}F(z) = F_D(z)
$$
and the product  converges when  $|z| < 1$. In addition,   $F(z) = \zeta_C(s)$ for $z=q^{-s}$ and the  function $F(z)$ will the function $Z_C(t)$ from the Grothendieck cohomology theory\footnote{A brief review of the theory see below, in section 4.1.} when we apply the substitution  $t=q^{-s}$.

From section 2, we know that
\begin{equation}\label{zeta}
\check\pi^*({\cal L}f_D) = z^{\deg D}\frac{P(z)}{(1-z)(1-qz)},
\end{equation}
and we may take as the target space in ($\ref{eq-emb}$) the space $\mathbb C_{++}[\mathbb T_0]$  of rational functions with possible first order poles in
the points $z = 1$ and $z = q^{-1}$.
Moreover, the image of the map $\check\pi^*$  can be identified with the space  $P(z)\mathbb C_{++}[\mathbb T_0]$. This is easy to check by applying
the map $\check\pi^*$ to the standard functions.

The first main result of these notes is the following statement:
the diagram
\begin{equation}\label{eq-comm}
\begin{CD}
\C_+[\mathbb T] @> F >>  \C_+[\mathbb T]\\
@ V \check\pi^*VV   @ V \check\pi^*VV\\
\mathbb \C_{++} [\mathbb T_0] @ > i^* >> \mathbb \C_{++}[\mathbb T_0]
\end{CD}
\end{equation}
commutes. Here $i: \mathbb T_0 \rightarrow \mathbb T_0$ is the involution $i(z) = q^{-1}z^{-1}$ of the torus $\mathbb T_0$.
To obtain  this commutativity we can consider the maps for the standard functions and verify the property for them. Namely, we then have a
diagram
%$$
%\xymatrix@!0{
%K_P\ar@{-}[dd] & & &  & K_{P,C}\ar@{-}[dl] \ar@{-}[dr]  &
%\\
% & & & K_P\ar@{-}[dr]  & & K_C \ar@{-}[dl]
%\\
%K & & & & K &
%}
%$$

%$$\xymatrix@{
%f_D = \underset x\otimes  f_{x,D}= \underset x\otimes z_x^{-v_x(D)}(1-z_x)^{-1}   \ar@{ \mid - >}[d]^{\Fo} \ar@{ \mid - >}[r]^{j^*} & F_D = %z^{\deg D} Z_C(z) \ar@{ \mid - >}[d]^{i^*}\\
%\underset x\otimes q^{-\frac 12 k_x+v_x(D)}f_{x, {\rm K}-D} \ar@{ \mid - >}[r]^{j^*} &q^{-\frac 12\deg {\rm K}+\deg D}z^{\deg {\rm K}-\deg D} %Z_C(z)}
%$$

$$
\begin{CD}
f_D = \underset x\otimes  f_{x,D}= \underset x\otimes z_x^{v_x(D)}(1-z_x)^{-1} @> \check\pi^* >> F_D = z^{\deg D} Z_C(z)\\
@ V \Fo VV  @V i^*VV\\
\underset x\otimes q_x^{-\frac 12 k_x+v_x(D)}f_{x, {\rm K}-D}
 @> \check\pi^* >> q_x^{-\frac 12\deg {\rm K}+\deg D}z^{\deg {\rm K}-\deg D} Z_C(z)
\end{CD}
$$
Since $i^* F_D(z)=F_D(q^{-1}z^{-1})$
we have to check the following claim:
$$
F_D(q^{-1}z^{-1})=q^{-\frac 12\deg {\rm K}+\deg D}F_{{\rm K}-D}(z).
$$
Indeed,

\begin{align*}
&F_D(q^{-1}z^{-1})=q^{\deg D}z^{\deg D}P(q^{-1}z^{-1})(1-q^{-1}z^{-1})^{-1}(1-z^{-1}) =\\
& q^{\deg D} z^{\deg D}[q^{-g} z^{-2g}] P(z)] qz^2 (1-qz)^{-1}(1-z)^{-1} =\\
&q^{-\frac 12\deg {\rm K}+\deg D}z^{-\deg {\rm K}+\deg D}P(z)(1-z)^{-1}(1-qz)^{-1} =\\
&q^{-\frac 12\deg {\rm K}+\deg D}F_{{\rm K}-D}(z)
\end{align*}
where $g = g(C)$ and we have used equality (\ref{zeta}) and a relation for  $i^*P$ following from the functional equation
\begin{equation}\label{eq-fu}
\zeta_C(s) = q^{g(C)-1}q^{-s(2g-2)}\zeta_C(1-s).
\end{equation}
for $\zeta_C(s) = Z_C(q^{-s})$ with $z = q^{-s}$.

In other words,  in each local space the Fourier transform has  additional multipliers which contribute  to keeping  the singularities of the functions
at the prescribed  points $1$ and $q_x^{-1}$. The local Fourier transform does not correspond to any geometrical automorphism of the
torus ${\mathbb T}_x$.  The situation will drastically change  when we go to the global space. In this case,
 the additional multipliers which we have for each local Fourier transform will disappear when we multiply
them over all the points of our curve $C$. Indeed, we have
$$\left.\prod_x \frac{1-q_x^{-1}z_x^{-1}}{1-z_x}\right|_{z_x=z^{\deg x},~ q_x=q^{\deg x}} =
Z_C(z) Z_C(q^{-1}z^{-1})^{-1} = q^{\frac12\deg {\rm K}} z^{\deg {\rm K}} $$
and the last equality  is  exactly the functional equation.

We see that the functional equation is the symmetry with respect to the involution $i: z\mapsto q^{-1}z^{-1}$ once we
identify the coordinate $z$ on the torus $\mathbb T_0$ with the parameter $t = q^{-s}$ from the Grothendieck theory.
Then $i(z)=q^{-1}z^{-1}$ corresponds to the map  $s\mapsto 1-s$ since $z=q^{-s}$.

The last remark which we will make in this section  concerns  the function spaces and Fourier transforms on the connected component ${\mathbb T} _0$ of the torus $\mathbb T_{\Pic}$. The Mellin  transform connects the functions on the group $\Gamma_C/\Gamma_0 = \Z$  with  the functions on this torus. We have
the following  commutative diagram:
$$
\begin{CD}
 {\cal D}_{++}({\Z}) @> \rm M >> \C_{++}[\mathbb T_0]\\
@ V \Fo VV  @V \Fo VV\\
{\cal D}_{++}({\Z}) @> \rm M>> \C_{++}[\mathbb T_0]
\end{CD}
$$
where
$$
{\cal D}_{++}({\Z}) =\left \{f(n), ~n \in \Z :
\begin{array}{rccllcc}
&f(n)& =& 0&~ \mbox{for}~n& \ll & 0\\
\exists a, b \in \C~ \mbox{s. t.}~&f(n)& = &aq^n + b   &~ \mbox{for}~n &\gg& 0
\end{array}
\right\}.
$$
and
\begin{equation}\label{eq-fur}
(\Fo f)(n) = q^nf(-n), \quad\mbox{for}~n \in \Z,
\end{equation}

By our construction, this space is exactly the image $\pi_*({\cal D}_+(\Gamma_C)$, where $\pi : \Gamma_C \rightarrow \Gamma_C/\Gamma_0$ is the natural projection. We  have then the commutative diagram
\begin{equation}
\begin{CD}
{\cal D}_+(\Gamma_C) @> \pi_{*} >>  {\cal D}_{++}(\Gamma_C/\Gamma_0)\\
@ V {\rm M}VV   @ V {\rm  M}VV\\
\mathbb \C_{+} [{\mathbb T}_{\Gamma}] @ >{\check\pi}^* >> \mathbb \C_{++}[\mathbb T_0]
\end{CD}
\end{equation}

\bigskip

Up to now, we considered what happens with the functions on torus ${\mathbb T}_{\Gamma}$ when they restricted to the connected component
$\mathbb T_0$ of unity of the torus $\mathbb T_{\Pic}$. We have
$$
\mathbb T_{\Pic} = \prod_{\chi \in \check\Phi} \chi\mathbb T_{0}.
$$
The torus $\mathbb T_{\Pic}$ has the involution $i: \mathbb T_{\Pic} \rightarrow \mathbb T_{\Pic}$ which sends a character $\chi : \Gamma/\Gamma_l \rightarrow \C^*$ to the character $\breve{\chi}:=|\cdot|\chi^{-1}$ with  $|a| = q^{-\deg(a)}$. This involution coincides with the previous involution $i$ on the torus $\mathbb T_{0}$. We claim that again there exists a commutative diagram
\begin{equation}\label{eq-comm1}
\begin{CD}
\C_+[\mathbb T] @> F >>  \C_+[\mathbb T]\\
@ V \check\pi^*VV   @ V \check\pi^*VV\\
\mathbb \C_{++} [\mathbb T_{\Pic}] @ > i^* >> \mathbb \C_{++}[\mathbb T_{\Pic}].
\end{CD}
\end{equation}
By the definition, the map $\check\pi : \mathbb T_{\Pic} \rightarrow \mathbb T_{\Gamma}$  is given by the following formula
$$
\check\pi (\chi) = (z_x, x \in C)~\mbox{with}~z_x = \chi(t_x),
$$
where $\chi \in \mathbb T_{\Pic} \subset \Hom(\A^*, \C^*)$ and $t_x$ is a local parameter at the point $x \in C$ which defines an element $(\dots, 1, t_x, 1, \dots)$ in $\A^*$. Then, formally, the inverse image of the standard function
$\delta_{\geq 0}$  is equal to
$$
\prod_{x \in C }(1 - \chi(t_x))^{-1}.
$$
For any $\chi \in \mathbb T_{\Pic}$ we have the $L$-function $L_C(s, \chi)$ which is defined as the product
$$
L_C(s, \chi) = \prod_{x \in C }(1 - \chi(t_x)q_x^{-s})^{-1} =  \prod_{x \in C }(1 - \chi(t_x)t^{\deg x})^{-1},
$$
where $t = q^{-s}$. If $\chi|_{\Phi} \neq 1$, then the product converges for all $t$ and is a polynomial $P(t)$ in $t$ of degree
$2g$.  It satisfies  the functional equation
$$
P(t, \chi)  = \chi((\omega))q^{g-1}t^{2g-2}P(q^{-1}t^{-1}, \chi^{-1}),
$$
where $\omega \in \Omega^1_K, \omega \neq 0$ is a rational differential form on $C$ and $(\omega)$ is the canonical class in $\Pic(C)$.
This follows from the Tate-Iwasawa method since this $L$-function is a particular case of the general $L$-functions from section 2.3. Since our character
$\chi$ is an unramified one, we may set  $f = \delta_{\geq 0}$ in  (\ref{eq-int}). See also, Theorem 6 in \cite{W2}[ch. VII, n 7].

We easily obtain that $P(q^{-1}, \chi^{-1} ) = P(1, |\cdot|\chi^{-1})$ and thus $P(1, \chi)  = \chi((\omega))q^{g-1}P(1, |\cdot|\chi^{-1})$.
By definition, $\chi((\omega)) = \prod_x \chi(t_x)^{\nu_x(\omega)} = \prod_x z_x^{k_x}$ and $q^{g-1} = \prod q_x^{\frac{1}{2}k_x}$.
If we take into account that $P(1, \chi) = \check\pi (\delta_{\geq 0})|_{\chi\mathbb T_0}$, we get for this inverse image the following identity
$$
\check\pi (\delta_{\geq 0})|_{\chi\mathbb T_0}  = \prod_x z_x^{k_x} q_x^{\frac{1}{2}k_x} \check\pi (\delta_{\geq 0})|_{\breve\chi{\mathbb T}_0}.
$$
The same relation holds for the inverse images of the standard functions $f_D$ for any divisor $D$.  We can now repeat all the computations which were made above for the component $\mathbb T_0$ and get the diagram (\ref{eq-comm1}).

\bigskip

Let us  finally say a few words about the Weil problem for the functional $\Delta_\chi$ in this new setting.
The local component of the integral (\ref{eq-int})   will be an integral over $K^*_x$ that equals to the arithmetical progression for $(1-z_x)^{-1}$ and  converges for $|z_x| < 1$. Then the functional $\Delta_\chi(f_x)$ is well defined for all $z_x \neq 1$ by the evaluation formula
$$
\Delta_\chi\colon\mathbb C_+[{\mathbb T}_x]\to\mathbb C,\quad \Delta_\chi(f_x) =  f_x(z_x)
$$
where the character $\chi$ corresponds to the complex number $z_x \in {\mathbb T}_x$.
This will provide a solution to the Weil problem for the functions on the local tori. Indeed, the structure of ${\cal D}(\Z)$-module on the space ${\cal D}_+(\Z)$, defined by convolution (see section 2.3), goes to the structure of  ${\mathbb C}[{\mathbb T}_x]$-module on ${\mathbb C}_+[{\mathbb T}_x]$, defined by multiplication. On  $\mathbb C[{\mathbb T}_x]$ we have the functional  $\Delta_\chi$, defined by the evaluation formula for all $\chi$, and we are looking for $\Delta_\chi$ on ${\mathbb C}_+[{\mathbb T}_x]$ such that
$$
\Delta_\chi (f\cdot g) = f(z_x)  \Delta_\chi (g), ~\mbox{for all}~f \in {\mathbb C}[{\mathbb T}_x]~\mbox{and }~  g \in  \mathbb C_+[{\mathbb T}_x],
$$
see (\ref{weil1}). The solution for $\chi \neq 1$ was given above. When $\chi = 1$ and thus $z_x = 1$, we set
$$
\Delta_\chi(f_x) = \res_{(1)}( f(z_x)\omega),
$$
where $\omega = dz_x/z_x$. Concerning the  global case see section 3.5 below.
\subsection{The Poisson formula and residues}
The involution $i:\mathbb T_0 \rightarrow \mathbb T_0$ extends uniquely to a map on the projective compactification
$\overline{\mathbb T}_0 = \mathbb P^1$ of the torus $\mathbb T_0  = \mathbb C^*$.
Let  $\omega = dz/z$. Then
$$
i^*(\omega) = d(q^{-1}z^{-1})/q^{-1}z^{-1} =  dz^{-1}/z^{-1} = - \omega.
$$
Let $g \in \C_{++}[\mathbb T_0]$. The divisor of poles has the support
$$
\Supp (g)_{\infty} \subset \{z = 0\} \cup  \{z = q^{-1}\}  \cup \{z = 1\} \cup \{z = \infty\}.
$$
Thus,
$$
\res_{(0)}(g\omega) + \res_{(q^{-1})}(g\omega) + \res_{(1)}(g\omega) + \res_{(\infty)}(g\omega)  = 0.
$$
Since  $\res_P(\eta) =  \res_{i(P)}(i^*\eta)$ for any differential form $\eta$  and any point $P$,
\begin{equation}\label{eq-res-sum}
\res_{(0)}(g\omega) +  \res_{(1)}(g\omega) =   \res_{(0)}(i^*g\cdot\omega) + \res_{(1)}(i^*g\cdot\omega).
\end{equation}
The torus $\mathbb T_0$ is  a connected component of $\mathbb T_{\Pic}$ and we have
$$
\mathbb T_{\Pic} = \coprod_{\chi \in \check\Phi}\chi\mathbb T_0 \quad \subset \quad {\bar{\mathbb T}}_{\Pic},
$$
where the boundary of compactification ${\bar{\mathbb T}}_{\Pic}$ has the following structure
$$
\partial {\bar{\mathbb T}}_{\Pic} = \coprod_{\chi \in \check\Phi}\{ 0_{\chi}\} \coprod_{\chi \in \check\Phi}\{ \infty_{\chi}\}.
$$
Each component $\chi\mathbb T_0$ is isomorphic to $\C^*$ with a coordinate $z$ and the points $0_{\chi}$ and   $\infty_{\chi}$ are defined as $z = 0$ and $z = \infty$ respectfully.
The form $\omega$ on $\mathbb T_0$ is invariant under translations and thus defines a canonical invariant form
$\omega$ on $\mathbb T_{\Pic}$.

At last, we introduce the function space $\C_{++}[\mathbb T_{\Pic}]$ on  $\mathbb T_{\Pic}$ as the direct sum of the space $\C_{++}[\mathbb T_0]$ and the spaces  $\C[\chi\mathbb T_0]$ where $\chi \in \check\Phi$ and $\chi \neq 1$.

\bigskip

 We now can obtain our second main result. Let $f \in {\cal D}_+(\Gamma)$.  We then have a sequence of transformations
$$
f \mapsto  {\rm M}f \in \C_+[\mathbb T]_{\Gamma}  \mapsto  \tilde f := \check\pi^*{\rm M}f \in \C_{++}[\mathbb T_{\Pic}],
$$
where  the map $\pi$ is again the natural projection onto the larger group $\Gamma/\Gamma_l$ (see section 3.3).

We claim that  the equality (where the residues at the points $z = 1$ are performed on the torus $\mathbb T_0$)
$$
\sum_{\chi \in \check\Phi}\res_{(0_{\chi})}(\tilde f\omega) +  \res_{(1)}(\tilde f\omega) =   \sum_{\chi \in \check\Phi}\res_{(0_{chi})}(i^*\tilde f\cdot\omega) + \res_{(1)}(i^*\tilde f\cdot\omega)
$$
is  equivalent to the Poisson formula \eqref{eq-pf1}
for the function $f$. Since $\#\check\Phi = \#\Phi$, this follows from the  facts:
\begin{enumerate}
\item $i^*(\tilde f) = \widetilde{\Fo f}.$
\item $\sum_{\chi \in \check\Phi}\res_{(0_{\chi})}(\tilde f\omega) = \#\check\Phi~\delta_{\Gamma_l}(f).$
\item $\res_{(1)}(\tilde f\omega) = \#\Phi /(q-1)~\delta_{(0)}(f).$
\end{enumerate}

 The first property follows from the commutative diagrams (\ref{eq-comm}) and (\ref{eq-comm1})  we have considered above.

The second property.  For $f \in {\cal D}_+(\Gamma)$ we have $\delta_{\Gamma_l}(f) = (\pi_*f)(0)$. We
want to go from the  group $\Gamma/\Gamma_l$ to the dual torus $\mathbb T_{\Pic}$.

In general, for a discrete group $G$ and a function  $\varphi$ on  $G$ we have  $({\rm M}\varphi ) (\chi) = \sum_{g \in G} \varphi (g)\chi (g)$,
where $\chi \in \mathbb T$ and  $\mathbb T$ is the dual group.
Then,  the following equality holds
$$
\varphi (0) = \int_{\mathbb T}({\rm M}\varphi),
$$
where the integral can be defined as
$$
\int_{\mathbb T}\psi = (\#\mathbb T)^{-1}\sum_{\chi \in \mathbb T}\psi(\chi),
$$
when the group $\mathbb T$ is finite and
$$
\int_{\mathbb T}\psi = \res_{(0)}(\psi dz/z),
$$
when $\mathbb T \cong \C^*$\footnote{Compare with the maps $\alpha^*$ and $\beta_*$  in diagram (\ref{d}).}.  This immediately gives the property we need.

The third property. Let $f  = f_D = \otimes_{x \in C} \delta_{(\ge n_x)} =  \otimes_{x \in C}f_x$ for a divisor $D = \sum_x n_x\cdot x$.
Then  $\delta_{(0)}(f)  =  \prod_x \delta_{(0)}(f_x) = \prod_x\lim_{n \to +\infty}f_x(n)  = 1$. Also as we saw, $\tilde f =  z^{\mbox{deg}~ D}Z(z)$ and from the computation of the residue of the zeta-function at the point $s = 0$ i.e. $z = 1$ (see above, section 2. 5) we get that $\res_{(1)}(\tilde f\omega)
=  \#\Pic^0(C)(\Fb_q)/(q-1)$. The general case will follow by linearity, since  any function  $f$ is a finite sum $\sum_D c_D f_{D}$.

\bigskip

{\bf Scholium 2}. In the \'etale cohomology theory we have $\zeta_C(s) = Z_C(t)$ with $t = q^{-s}$ and the parameter
$t$ appears in a purely formal manner (see section 4.1 below). Our torus coordinate $z$ does appear in a canonical way as a coordinate on the
connected component of the torus dual to the group $\Gamma/\Gamma_0 = \Pic(C)$. We can say that $t = z$ and thus
the  function  $Z_C(t)$  is naturally defined on $\mathbb P^1 \supset \mathbb C^*$, not only on the torus $\mathbb C^*$.

\bigskip

We can also rewrite the more general Poisson formula in which  the shifted $\delta$-function appears.
Given  $g \in \A^*$,  then the analog  of the Poisson formula (\ref{eq-pf5}) with $\delta_{gK}$ will be the following
relation for residues:
$$
\res_{(0)}(z^n\tilde f\omega) +  \res_{(1)}(\tilde f\omega) =   \res_{(0)}(z^{-n}i^*\tilde f\cdot\omega) + \res_{(1)}(i^*\tilde f\cdot\omega)
$$
where $n = \deg{\gamma}$ and  $\gamma$ is the image of $g$ in the group $\Gamma_C = \A^*/\OO^*$.
%\subsection{The Weil problem on holomorphic tori (to be written)}
%The last thing that we have to do is to develop the whole Tate-Iwasawa method in the holomorphic setting. We saw above how to deal with the local
%Weil problem.
 \subsection{Explicit formulas}
The procedure we have carried out  here is very similar to the way in which  the so called explicit formulas are deduced in  analytic number theory.
We will present here the corresponding reasoning in  our language and  make a comparison  with traditional exposition in the next section.

Let  $P(z)(1 - z)^{-1}(1 - qz)^{-1} = Z(z) = \zeta_C(s)$  for $z = q^{-s}$. Then   $\omega = dZ/Z$ is a logarithmic differential form on the projective line $\P^1$ and
$$
\Supp (\omega)_{\infty} \subset \{z = 0\} \cup    \{z = q^{-1}\}  \cup  \cup_{\lambda} \{z = \lambda^{-1}\}  \cup  \{z = 1\} \cup   \{z = \infty\},
$$
where $P(z)  = \prod_{\lambda} (1 - \lambda z)$ is the spectral decomposition for the characteristic polynomial of the Frobenius automorphism.
Let $h$ belong to $\C[\mathbb T_0]$ and come from a function $f \in {\cal D}(\Z)$ via the Mellin  transform (see (\ref{l-tr})). Taking as above the sum of residues, we obtain
$$
\res_{(0)}(h\omega) + \res_{(q^{-1})}(h\omega) +  \sum_{\lambda}\res_{(\lambda^{-1})}(h\omega) + \res_{(1)}(h\omega)
 + \res_{(\infty)}(h\omega)  = 0.
$$
We want to compute all the members  of this sum:
\begin{enumerate}
\item $\res_{(q^{-1})}(h\omega) = - h(q^{-1}).$
\item $\res_{(\lambda^{-1})}(h\omega) =  h(\lambda^{-1}).$
\item $\res_{(1)}(h\omega) = - h(1).$
\item $\res_{(\infty)}(h\omega) =  \res_{(0)}(\widetilde{\Fo f}\omega) + (2 - 2g(C))\res_{(0)}(\widetilde{\Fo f} dz/z). $
\item $\res_{(0)}(h\omega) = - \sum_{x \in C, n \ge 1}(\deg x)f(-n\deg (x)).$
\end{enumerate}
The first three properties follow directly from computation of logarithmic residues and our knowledge of behavior of the function $Z(z)$ at the points
$z = q^{-1}, z = \lambda^{-1}, z = 1$.

To get the fourth  property we apply the involution $i$, which  transforms the point $z = \infty$ into the point $z = 0$, and use the invariance of the residue.
Next, the functional equation \eqref{eq-fu} implies that $i^*(\omega) = (2 - 2g(C))dz/z + \omega$ and our diagram \eqref{eq-comm} gives us that $i^*h =
\widetilde{\Fo f}$ if the original function $h$ was produced from the function $f$.

To get the fifth property we have to start with the Euler product $Z(z) = \prod_x(1 - z^{\deg x})^{-1}$.
We have in a neighborhood  of the point $z = 0$
$$
\omega = d\ln Z = - \sum_{x \in C, n \ge 1}(\deg x) z^{(n\deg x)}dz/z.
$$
If $h(z) = \sum a_mz^m$ then we arrive at the finite sum
$$
\res_{(0)}(h\omega) =  - \sum_{x \in C, n \ge 1}(\deg x) a_{-n\deg x},
$$
where $n\deg x \le -\nu_{(0)}(h)$.

Since we have started  from a function $f$ defined on the group $\Gamma/\Gamma_0 = \Z$, we have $a_m = f(m),  m \in \Gamma/\Gamma_0$ and
the final formula  will be
\begin{equation}\label{eq-weil-explicit}
\begin{array} {l}
h(q^{-1}) -  \sum_{P(\lambda^{-1}) = 0} h(\lambda^{-1}) +   h(1) =\\[4mm]
(2 - 2g(C))\Fo f(0) - \sum_{x \in C, n \ge 1}(\deg x)[ f(-n\deg x)  + \Fo f(-n\deg x)]   =\\[4mm]
(2 - 2g(C)) f(0) - \sum_{x \in C, n \ge 1}(\deg x)[ f(-n\deg x)  + q^{-n\deg x} f(n\deg x)] .
\end{array}
\end{equation}
Here, we could replace  $\Fo f(0)$   by $ f(0)$ and $\Fo f(n)$  by $q^nf(-n)$ .

This is exactly the explicit formula for the case of function fields (see, \cite{W3, C} and compare with the case of number fields in \cite{L1, VK}).
The first explicit formula was found by Riemann. It relates the distribution of prime numbers to a sum over critical zeros of the zeta function. Later de la Vall\'ee Poussin and Hadamard   added some arguments from complex analysis and, after that,  this formula leads to the  prime number theorem. In our case, the closed points $x$ replace the prime numbers and the eigenvalues  $\lambda$ replace the critical zeros of zeta function. The corresponding asymptotic law is also valid and can be stated  as follows:
$$
 \# \{ x \in C: N(x) < N\}  \sim \frac{N}{\ln N}~\mbox{as}~N \to \infty,
$$
where $N(x) = q_x = q^{\deg x} = \# k(x)$.

To obtain this  from the explicit formula (26) one chooses for $f$ the test function $f = \delta_{\ge N}$. The  de la Vall\'ee Poussin and Hadamard  reasoning was to show that there are no zeros of the zeta-function on the boundary $\Re(s) = 1$ of the critical strip. In our case,  it is equivalent to the inequality $|\lambda| < q$
(the critical strip will be $q^{-1} \le z \le 1 $) which is much weaker then the Riemann's conjecture $|\lambda| = q^{-1/2}$ (Weil's theorem) .

Let us also note that by a simple computation the asymptotic law is equivalent to the following asymptotic behavior  for the number of rational points on the curve $X$:
$$
\frac{X({\mathbb F}_{q^n})}{q^n} \to 1\quad \mbox{when}~n \to \infty.
$$

\section{Trace formula and the Artin representation}
\subsection{Zeta functions and \'etale cohomology}
Let $X$ be a scheme of finite type over $\Z$. We denote by $|X|$ the set of closed points of $X$.
For every $x\in |X|$ the residue
field $k(x)$ is   a finite field.
The zeta function of  $X$ is
defined  as
$$
 \zeta(X,s)\ := \prod_{x\in |X|}(1- (\#k(x))^{-s})^{-1} \ , s\in\C \  .
 $$
This expression is called the Euler product for the zeta-function.
The Euler product
converges for ${\rm Re}(s)>\dim(X)$.

The first examples are the following:
\begin{itemize}
 \item[1.] $\zeta(\Spec(\Fb_q),s)= (1-q^{-s})^{-1}$.
 \item[2.] $\zeta(\Spec(\Z),s)=\zeta(s)$,  the Riemann zeta function.
 \item[3.] Let $X/\Fb_q$ be a scheme over a finite field. Then
           \begin{equation}\label{1f}
            \zeta(X,s)\ =\
               \exp \left( \sum_{m=1}^\infty \# X( {\mathbb F}_{q^m} ) \frac{q^{-ms}}{m}
                    \right)
           \end{equation}
           in the ring  of formal power series and the series is absolutely  convergent
           for ${\rm Re}(s)>\dim(X)$.
\end{itemize}
To prove 3.  one takes the logarithmic derivative of both sides in the ring of formal
power series and uses the relation
\begin{equation}\label{numbers}
\# X({\mathbb F}_{q^m})\ =\ \sum_{l | m} a_l\, l \ ,\ m=1,2,\ldots \  ,
\end{equation}
where $a_l$ denotes the number of closed points $x\in X$ for which $\#k(x) =
q^l$. This number is finite. Also, the absolute convergence of the infinite
product and of the power series are equivalent.

If now $X={\mathbb A}^n_{\Fb_q}$, we see that
$ \zeta(X,s)\ =\ \exp \left( \sum_{m=1}^\infty q^{m(n-s)}/m \right) $
and the sum on the right-hand side is absolutely convergent for $|q^{n-s}|<1$,
i.e. for ${\rm Re}(s)>n=\dim(X)$.
Furthermore, we see that
\begin{itemize}
 \item[4.] $\zeta({\mathbb A}^n_{\Fb_q}, s)\ = \ (1-q^{n-s})^{-1} \  .$
\end{itemize}
In the general case, the convergence of the Euler product can  be easily reduced to this case.
For projective varieties over $\Fb_q$, the main tool for the  study of the zeta (and also $L$-) functions
 is the  \'{e}tale cohomology, or more precisely  $l$-adic cohomology for some prime number $l$ such that  $(l, q) = 1$.
Let us fix such an $l$ and a ground field $k$.

For every projective variety (actually, for much more general schemes)  $X/k$,  there exists  a family
$$ H^i(X,\Ql),\ i=0,1,\ldots $$
of $\Ql$-vector spaces called the $i$-th $l$-adic cohomology groups  with values in $\Ql$.
The groups $H^i(X,\Ql)$ are finite dimensional vector spaces over $\Ql$
and  are trivial for $i> 2 \dim(X)$ when the ground field $k$ is algebraically closed.

The cohomology  group $H^i(X,\Ql)$ is contravariant with respect to morphisms of varieties
and there exist the usual cup product pairings which are non-degenerate (see
[{\bf SGA4, SGA4-1/2, SGA5}] for the original sources, as well as the standard  textbooks).

When $k=\C$,  these groups coincide with the groups $H^i(X,\Q)\otimes_\Q \Ql$, where
$H^i(X,\Q)$ are the usual cohomology groups of the topological space $X(\C)$ (for the classical topology) and  many
properties of the topological cohomology groups remain true in the $l$-adic setting.
In particular, there is the Lefschetz fixed point formula for a morphism $f: X \rightarrow X$
of projective varieties defined over an algebraically closed field.  We have
\begin{equation}\label{eq-lef}
\sum_i(-1)^i\Tr (f: H^i(X,\Ql) \rightarrow H^i(X,\Ql))  =  \sum_{x \in X}  \# (\Gamma_f \cap\Delta)_{(x,x)},
\end{equation}
assuming that  the graph $\Gamma_f$
of the morphism  and  the diagonal $\Delta \subset X \times X $ are in general position in $X \times X$. If we  assume that  the subvariety $\Gamma_f$ intersects $\Delta$  transversally, then
all the multiplicities at the fixed points $x \in X$ of $f$  equal to 1 and we can conclude that the formula counts the number of  the fixed points:
\begin{equation}\label{eq-lef-1}
\sum_i(-1)^i\Tr (f: H^i(X,\Ql) \rightarrow H^i(X,\Ql))  = \sum_{x \in X: f(x)=x}1.
\end{equation}
In order to show how this formula can be applied to the study of zeta-functions
one needs the Frobenius automorphism.
It is well defined as an automorphism $Fr:\bar{X}\to\bar{X}$ of the scheme  $\bar{X} =  X \otimes_{\Fb_q}\bar{\Fb_q}$ over  $\bar{\Fb}_q$ for
any  scheme $X$ over a finite field $\Fb_q$.
For the scheme $X$ over $\Fb_q$, we define $Fr_X$ as identity on the topological space  $X$ and $q$th power on the structure sheaf ${\cal O}_X$.
On the scheme $\bar{X}$, the Frobenius  morphism $Fr$ is $Fr_X \otimes 1_{\bar{\Fb_q}}$. If the scheme $X$ is embedded into a projective space
with coordinates  $(x_0 : x_1 : \ldots : x_n)$ then we have an embedding of the scheme $\bar{X}$ into the same projective space. The  Frobenius morphism will be then a restriction to     $\bar{X}$ of the map  $(x_0 : x_1 : \ldots : x_n) \mapsto (x_0^q : x_1^q : \ldots : x_n^q)$.
The following properties hold :
\begin{itemize}
 \item[1.] $\bar{X}^{Fr^n} = X({\mathbb F}_{q^n})$, where $\bar{X}^{Fr^n}$ denotes the
           subset of closed points of $\bar{X}$ which are fixed by  $Fr^n$ .
 \item[2.] The set of closed points $x$ of $X$ can be identified with the set of
           orbits of $Fr$ acting on the set of closed points of $\bar{X}$.
           The integer $\deg(x) = \dim_{\Fb_q} k(x)$ is equal to the number of elements
           in the orbit corresponding to $x$.
 \item[3.] $\# \bar{X}^{Fr^n} = \sum_{x \in X, \#k(x)|n} \deg(x)$.
\end{itemize}

 In our case of projective varieties, we have
an induced action
$$ Fr^*\ :\ H^i(\bar{X},\Ql)\ \to\ H^i(\bar{X},\Ql) \ . $$
The differential of $Fr$  is  the  zero  map. This implies that the graph of $Fr$ intersects the diagonal
transversally, and the Lefschetz fixed point formula for  $l$-adic
cohomology now states that for $n = 1, 2, \dots$
$$ \# {\bar X}^{Fr^n}\ =\
    \sum_{0 \le i \le2\dim(X)} (-1)^i \Tr ( (Fr^*)^n\,:\, H^i(\bar{X},\Ql)\,\to\,H^i(\bar{X},\Ql) ) \ .
$$
For an arbitrary linear operator $F$ acting on a finite dimensional vector space $V$
we have the following relation in the ring of formal power series:
$$ t \frac{d}{dt} \log \left( \det ( 1-F t ) \right)^{-1} \ =\
                        \sum_{n=1}^\infty \Tr(F^n) t^n \
\footnote{For the proof one can assume the base field to be algebraically closed.
 One then checks the $1$-dimensional case directly and uses additivity of
the both sides of the formula with respect to short exact sequences of vector spaces $V$.}.
$$
Using the relation \eqref{1f} it makes sense to set $\zeta(X,s)= Z(X,t)$
with $t=q^{-s}$. We get
%$$\begin{array}{rcl}
  % t \frac{d}{dt} \log  Z(X,t)
   %&=& \sum_{n=1}^\infty \# X({\mathbb F}_{q^n}) t^n                            \\[1mm]
   %&=& \sum_{i=0}^{2\dim(X)} (-1)^i
     %         \sum_{n=1}^\infty \Tr( (F^*)^n:
       %       H^i(\bar{X},\Ql)\ \to\ H^(\bar{X},\Ql) )   \\[1mm]
   %&=& \sum_{i=0}^{2\dim(X)} (-1)^i t \frac{d}{dt} \log \left(
      %              \left( \det ( 1-F^* t , H^i(\bar{X},\Ql) ) \right)^{-1} \right) \ .
  %\end{array}
%$$
$$
t \frac{d}{dt} \log  Z(X,t)   = \sum_{n \ge 1} \# X({\mathbb F}_{q^n}) t^n =
$$
$$
 \sum_{i \ge 0} (-1)^i  \sum_{n \ge 1} \Tr( (Fr^*)^n\, : \,
              H^i(\bar{X},\Ql)\ \to\ H^i(\bar{X},\Ql) )   =
 $$
$$
 \sum_{i \ge 0} (-1)^i t \frac{d}{dt} \log \left(
                    \left( \det ( 1-Fr^* t \, \vert H^i(\bar{X},\Ql) ) \right)^{-1} \right)  .
$$
Therefore,
\begin{equation}\label{1g}
  Z(X,t)\ =\ \prod_{0 \le i \le 2\dim(X)} \left( \det ( 1-Fr^* t\, \vert H^i(\bar{X},\Ql) )
                                            \right)^{(-1)^{i+1}} = P(t)/Q(t) \ .
\end{equation}
This gives us a description of the zeta function in terms of $l$-adic cohomology and
we see immediately that $\zeta(X,s)$ is a rational function in $q^{-s}$
since  the  spaces $H^i(\bar{X},\Ql)$ are finite-dimensional.

In the case of an algebraic curve $C$, we have
$P(t) =  \det ( 1-Fr^*  t \, \vert H^1(\bar{C},\Ql) )$ (a polynomial of degree $g(C)$) and
$Q(t) = \det ( 1-Fr^*  t \, \vert H^0(\bar{C},\Ql) )\det ( 1-Fr^*  t \, \vert H^2(\bar{C},\Ql) ) = (1-t)(1-qt)$,
as we have already found above.

The last remark is that  Poincar\'{e} duality holds  in the $l$-adic
cohomology and implies the functional equation for zeta functions.

\subsection{Local Artin representation}
Let $K$ be a local field (of dimension 1), i.e. the field of fractions of a complete discrete valuation ring ${\cal O}$, and $L/K$  a finite normal
extension. Then $K \supset {\cal O}_K \supset {\mathfrak m}_K$ with the maximal ideal ${\mathfrak m}$ of  ${\cal O}_K$ and $k$ the residue field of $K$.  Let $l$ be the residue field of $L$. Let $g \in \Gal(L/K)$. We set:
$$
i_G(g) = \mbox{maximal integer}~ i ~\mbox{such that}~g~\mbox{acts trivially on}~{\cal O}_L/{\mathfrak m}_L^{i+1}.
$$
We have $i \ge -1$ and $i_G(e) = \infty$.The function $i_G$ defines a filtration $\{G_i\}$ of the group $\Gal(L/K)$. An extension is unramified if and only if
$G_0 = \{e\}$,  i.e. if $i_G(g) = -1$ for all $g \neq e$.

If $g \in G_0$ then we can redefine the function
$$
i_G(g) = \nu_L(g(t) - t), \quad g \neq e,
$$
where  $\nu_L$ is the discrete valuation of the field $L$ and $t$ is a generator of the maximal ideal  ${\mathfrak m}_L$.

The basic observation is that if we introduce the new function $a_G$ on $G$
$$
a_G(g) = \left\{
\begin{array}{rll}
-& i_G(g)[l:k] , & g \neq e,\\
\sum_{g \neq e}&i_G(g)[l:k], & g = e,
\end{array}
\right.
$$
then this function will be a character of some (finite-dimensional) representation of the group $G$ over $\C$. This fact can be proved in a very non-direct way  and it is still an open problem to define this representation, called the Artin representation,  in an explicit manner (see \cite{Se2}).

Let the residue field of $K$ be a finite field $\Fb_q$ and the field $K$  be a field of equal characteristics. There is class field theory describing the Galois group of the maximal abelian extension of $K$ \cite{A, W2, Se2}.
We have the commutative diagram
$$
\begin{CD}
{\cal O}^*   @ > >>K^*@ >\deg >> \Z\\
@ V\varphi  VV   @ V\varphi VV @ V\varphi VV\\
\Gal(K^{ab}/\bar\Fb_qK) @>>>  \Gal(K^{ab}/K)  @>>>\Gal(\bar\Fb_qK/K)
\end{CD}
$$
The system of groups  $G^i : = G_{\psi(i)}$  where $\psi : \Z \rightarrow \Z$ is the Herbrand  function (see, \cite{A, W3, Se2}) will be consistent in the tower of all finite normal extensions of the field $K$. Thus, we can  define the Artin  representations of  infinite extensions, including the maximal abelian extension. The reciprocity map $\varphi$ has the following property:
$$
\varphi (1 + {\mathfrak m}^i) =  \Gal(K^{ab}/K)_{\psi(i)}.
$$
 This means that the Artin character can also be  defined on the group $K^*$.  The corresponding function $a_K$ can then be considered as a distribution on the locally compact group  $K^*$ and  $a_K \in {\cal D}'(K^*)$. The definition is rather subtle but still very explicit:
$$
a_K(f)  = \int\limits_{{\cal O}^*} \frac{f(g) - f(1)}{|1 - g|} d^*g + \int\limits_{K^* - {\cal O}^*} \frac{f(g)}{|1 - g|} d^*g, \quad f \in {\cal D}(K^*).
$$
The integral without the substraction of a term $f(1)$  can be divergent if $f(1) \neq 0$  and the given formula contains a regularization  required in order to define the integral  for any function belonging to  the space ${\cal D}(K^*)$  (see \cite{W3, C}).

Assume that the function $f$ is invariant under ${\cal O}^*$. Then the first term will be zero, $K^*  = \coprod_nt^n{\cal O}^* $,  and for the full character we get
$$
a_K(f)  =  \int\limits_{K^* - {\cal O}^*} \frac{f(g)}{|1 - g|} d^*g =
$$
$$
\sum_{n \ge 1}\int_{t^n{\cal O}^*}f(t^n)|1 -g|^{-1}d^*g
+ \sum_{n \le -1}\int_{t^n{\cal O}^*}f(t^n)|1 - g|^{-1}d^*g =
$$
$$
\sum_{n \ge 1}f(t^n)\vol({\cal O}^*) + \sum_{n \le -1}\int_{t^n{\cal O}^*}f(t^n)q^{-n}\vol({\cal O}^*) =
\vol({\cal O}^*) \sum_{n \ge 1}[f(t^n) + f(t^{-n})q^{n}].
$$
Here, we have used that $|1 - g| = 1$ if $|g| < 1$ and $|1 - g| =  |g|\cdot |1 - g^{-1}| = |g|$ if $|g| > 1$.
\bigskip

Now we return to the global case. Let $X$ be an algebraic curve (as usual irreducible, projective and smooth) defined over a finite field  $\Fb_q$ and let $G$ be a finite group of  automorphisms of $X$. Also, we have the curve $C = X/G$ and a finite covering map $\varphi : X \rightarrow C$.   If $L = \Fb_q(X), K = \Fb_q(C)$ are the fields of rational functions on the curves then $G = \Gal(L/K)$ and for any point $y \in X$ we have the decomposition group $D_y \subset G$, defined by  $D_y = \{g \in G : g(y) = y \}$. If $x = \varphi (y)$ there is  a finite normal extension of  local fields $L_y/K_x$ with the Galois group $D_y$.
Thus, we have the functions $i_G$ and $a_G$  on the groups $D_y$ which we denote by  $i_y$ and $a_y$.

We can show that
$$
i_y(g) = \nu_y(g(t) - t) = (\Gamma_g\cdot\Delta)_{(y, y)}.
$$
 where $\Gamma_g$ is the  graph of the automorphism $g$ in the surface $X\times X$ and $\Delta$ is the diagonal. The first equality has already been mentioned. The second one follows from the definition of  the intersection multiplicity   of $\Gamma_g$ with $\Delta$ at the point $(y, y)$. Indeed, in local coordinates $u, t$  an equation for  $\Gamma_g$ is given by $f_1 = g(u) - t$  near the point $(y, y)$ and  an equation for  $\Delta$ by $f_2 = u - t$. By definition, $(\Gamma_g\cdot\Delta)_{(y, y)} = \dim {\cal O}_{(y,y)}/ (f_1, f_2)$ and $\dim {\cal O}_{(y,y)}/ (f_1, f_2) = \dim {\cal O}_y/(g(u) - u)$.

 For all $y$ lying over a given $x \in C$ the functions $a_y$ can be extended by zero to the whole group $G$ and then we can set
 $$
 a_x :  = \sum_{y \mapsto x} a_y.
 $$
For any point $y$ lying over $x$, the character  $a_x$ will be induced by  $a_y$. The representation $\pi_x$  of the Galois group $G$ defined by the character $a_x$ will be called the Artin representation at the point $x \in C$.

We arrive at  the finale of this construction, an application of the Lefschetz trace formula:
\begin{equation}\label{eq-lf2}
\Tr (g)|H^{\cdot}(\bar X, \Ql)   = (2 - 2g(C))\# G\delta_G(g) - \sum_{x \in C}a_x(g),
\end{equation}
where $\delta_G(g) = 0$ for $g \neq e$ and  $\delta_G(e) = 1$.

{\em Proof.} For $g \neq e$ this is exactly the Lefschetz trace formula \eqref{eq-lef} applied to the extension to the curve $\bar X$ of the automorphism $g$.  Also, use the computation of the multiplicities made above. If $g = e$ then we have the  formula  $(2g(X) - 2) = [L:K](2g(C) - 2) +$  (local terms) for the behavior of the genus under a finite covering and after  base change of the ground field \cite{T2}. For example, if our extension is only weakly ramified (that is  $G_i = \{e\}$ for $i > 1$) then we have the Hurwitz formula
$$
(2g(X) - 2) = [L:K](2g(C) - 2) + \sum_{x \in C} (e_x - 1),
$$
where $e_x$ is the ramification index over the point $x$.

Basically this result was found by Andre Weil.  Here, we have followed  Serre's exposition in \cite{Se2}.

\subsection{Relation with the explicit formulas}
We now want to combine  the explicit formula of section 3.4 and the trace considerations we have just discussed. Let us first review the explicit formulas, in  Weil's exposition \cite{W2}.  We consider again an algebraic curve $C$, the field $K = \Fb_q(C)$ and its maximal abelian extension $K^{ab}$.  By the class field theory, there exists the reciprocity map $\A_C^*/K^* \rightarrow \Gal(K^{ab}/K)$ and, also,  the degree map $\deg: \A_C^*/K^* \rightarrow \Z$. We start with a function $f = f(n), n \in \Z$, where the group $\Z$ is the target of the degree map. Let us assume that $f \in {\cal D}(\Z)$ so that  its Laurent transform $\tilde f$ belongs to $\C[\C^*]$. With Weil,  we introduce the functional
$$
S(f) =  \sum_{P(z) = 0}\tilde f(z),
$$
where $P$ is the numerator of the zeta function of $C$. The degree map allows us to consider the function $f$ as a function on the group $\A_C^*/K^*$.
Also, since we have embeddings $K^*_x \subset \A_C^*/K^*$ we can restrict the function  $f$  to $K^*_x$ and get a function  $f_x$ for each  point $x \in C$. We have the following distributions
$$
D_x(f) =  \int\limits_{{\cal O}_x^*} \frac{f_x(g) - f_x(1)}{|1 - g|} d^*g + \int\limits_{K_x^* - {\cal O}_x^*} \frac{f_x(g)}{|1 - g|} d^*g,
$$
$$
D(f) = \int\limits_{\A_C^*/K^*}f(g) [|g|^{1/2} + |g|^{-1/2}] d^*g.
$$
In our setting, the Weil explicit formula reads
\begin{equation}\label{eq-weil-expl}
S(f) = (2 - 2 g(C)) \delta_{e}(f) + D(f)  - \sum_{x \in C}D_x(f).
\end{equation}
By the definition,  we have $D_x(f)  = a_x(f_x)$ and  therefore  we can compare the expressions \eqref{eq-weil-expl}
and \eqref{eq-lf2}. For this , one needs to rewrite  \eqref{eq-lf2} for functions of the elements  $g \in G$ .We see that $S(f) - D(f)$ will correspond to the Lefschetz trace in the cohomology. The first term $S(f)$   corresponds to the trace on the group
$H^1(X, \Ql)$, and the second term  $D(f)$  to the traces on the groups  $H^0({\bar X}, \Ql)$ and $H^2({\bar X}, \Ql)$.

The works \cite{C, M, G} contains an interpretation of  all terms in the Weil formula as  traces. Connes uses representations in the Hilbert spaces, Meyer in the Bruhat-Schwartz spaces. Here, we give our own, much more elementary,  version.

Consider  the Lefschetz formula for integer powers of the Frobenius map. For any $n \ge 1$, we have (by (\ref{eq-lef-1})):
\begin{equation}\label{eq-lef-2}
\begin{array}{lllcl}
\# C(\Fb_{q^n})& = & \# {\bar C}^{Fr^n}& =&     \Tr ( (Fr^*)^n\vert~ H^{\bullet}(\bar{C},\Ql) \\
&&2 - 2g(C) & =&     \Tr ( (Fr^*)^0\vert~ H^{\bullet}(\bar{C},\Ql) \\
q^{-n}\# C(\Fb_{q^n})& = & q^{-n}\# {\bar C}^{Fr^n}& =&     \Tr ( (Fr^*)^{-n}\vert~ H^{\bullet}(\bar{C},\Ql) .
\end{array}
\end{equation}
Let $f: \Z \rightarrow \C$ be a function on $\Z$ with  finite support. Recall  that by the class field theory $1 \in \Z = \A^*/\A^{*(1)}$ corresponds to the Frobenius automorphism of the field $\bar{\Fb_q}$ over $\Fb_q$. We may thus introduce
$$
f(Fr) :  = \sum_{n\in \Z}f(n)Fr^{-n},
$$
 and conclude that
 $$
 \sum\nolimits_{n \ge 1}f(n)q^{-n} \# C(\Fb_{q^n})  + (2 - 2g(C)) f(0) + \sum\nolimits_{n \le -1}f(n)\# C(\Fb_{q^{-n}}) = \Tr ( f(Fr)\vert~ H^{\bullet}(\bar{C},\Ql) .
  $$
 As usual, we denote by $\tilde f (z) = \sum\nolimits_n f(n)z^n$ the Laurent transform of the
function $f$.
The RHS  of the equation is equal to
$$
\tilde f (1) + \tilde f (q) - \sum_{P(z) = 0}\tilde f(z).
$$
On the LHS, we get
$$
 \sum\nolimits_{n \ge 1}(\sum_{x \in C(\Fb_{q^n})}1 )q^{-n}f(n) +
(2 - 2g(C))f(0) + \sum\nolimits_{n \le -1}(\sum_{x \in C(\Fb_{q^{-n}})}1)f(n) =
$$
$$
(2 - 2g(C))f(0) + \sum\nolimits_{n \ge 1}\sum_{ x \in C, \deg(x)\vert n}\deg (x)q^{-n}f(n) + \sum\nolimits_{n \ge 1}\sum_{ x \in C, \deg(x)\vert n}\deg (x)f(-n) =
$$
$$
(2 - 2g(C))f(0) + \sum\nolimits_{n \ge 1, x \in C}\deg (x)q^{-n\deg(x)}f(n\deg (x)) + \sum\nolimits_{n \ge 1, x \in C}\deg (x) f-(n\deg (x)).
$$
Here, we used the fact that, according to  (\ref{numbers})  from section 4.1,  there is an equality $ \# {\bar C}^{Fr^n} = \# C(\bar{\Fb}_{q^n}) = \sum_{x \in C, \deg (x)|n}\deg (x)$.
We arrive once more at the same explicit formula (26) that was explained above.

\section{Number fields (from $\C^*$  to $\C$)}
It is certainly interesting to understand what corresponds to the constructions developed above in the case of  number fields. According to fundamental analogy between number fields (= finite extensions of $\Q$) and  fields of algebraic functions of one variable over a finite field (= fields of functions on algebraic curves) (see, for example, \cite{P2}) we expect that there will  be some analogous constructions. Actually, the adelic approach to the zeta- and $L$-functions developed by Tate and Iwasawa has realized this idea with respect to the analytic properties of these functions (see especially the unified exposition by Weil \cite{W1}\footnote{At the same time,  I. R. Shafarevich gave a course of lectures at  Moscow University where a unified exposition was also presented \cite{Sh}. An interesting historical coincidence!}). Thus what was said  in section 2 can be transfered to the number field situation with appropriate modifications\footnote{Historically, the path was in the opposite direction!}. The main problem is to deal with the archimedean valuations of the fields.

If we wish to reconsider our duality constructions for the number fields then we must first note that there is no the substitution  $t = q^{-s}$ for the zeta-function since
we no longer have  any $q$. The $L$- functions  are very simple analytic functions in the case of curves, since  they are superpositions of a rational function and the exponential function. This is already no longer     for the Riemann's zeta function. This  means that we have to reformulate our results on the $s$-plane, forgetting the variables such as $z$.

The map $s \mapsto z = q^{-s}$ is a universal covering of the torus $\C^*$ under an action of the discrete group $\Z$ and the $s$-plane can be considered as $\Hom(\Z, \C)$ if we start from the group $\Z$. The action is generated by the map $s \mapsto s + 2\pi i/\ln q$. We will discuss this issue in more detail later in the text and will now  simply rewrite our residue sum on the
$s$-plane.

Let $f = \delta_{(\ge 0)} \in {\cal D}_+(\Gamma)$ for the  group $\Gamma_C$ associated to  a curve $C$. Then $\tilde f =  Z_C(z)$. If we rewrite the sum of residues for $\tilde f\omega_0$
on the $s$-plane we meet the problem that there are infinitely many poles  lying on arithmetical progressions. The easiest way to overcome this difficulty
is to consider a fundamental domain for the group $\Z$. We may take a horizontal strip $ \pi/\ln q \ge t \ge -\pi/\ln q$. The image of this strip will be exactly
the open subset $\C^*$ of   $\mathbb P^1$. We need to add to the strip two infinite points $\pm\infty$ and then we must have
$$
\res_{(+\infty)}\zeta_C(s)ds + \res_{(0)}\zeta_C(s)ds + \res_{(1)}\zeta_C(s)ds + \res_{(-\infty)}\zeta_C(s)ds = 0.
$$
In the "open" strip, there are only two poles and we can observe  that the integral over the rectangle $\gamma$ with edges
$ \gamma_1,~ \gamma_2,~ \gamma_3,~ \gamma_4$
$$
\begin{array}{llllllll}
\gamma_1 & = & \{A + it & : &  -\pi/\ln q & \leq& t &\leq \pi/\ln q \},\\
\gamma_2 & = & \{\sigma + iT& : & -A &\leq &\sigma &\leq A \},\\
\gamma_3 & = & \{-A + it &: &  -\pi/\ln q &\leq &t &\leq \pi/\ln q \},\\
\gamma_4 & = & \{\sigma + iT& : & -A &\leq &\sigma &\leq A \},
\end{array}
$$
where $s=\sigma+it$ and $A$ is sufficiently large,
will be zero. Of course, the integrals over vertical boundaries of the strip do not depend on the choice of $A$. Thus the sum makes sense without any
additional infinite points and, certainly, the sum defined as this integral will be equal to our previous sum of residues on the $z$-torus.
 This suggests what to do in the case of a number field $K$. Denote by  $\xi_K(s)$ the 'complete'  zeta-function of the field $K$. This is the Dedekind's
 zeta function completed by the gamma-multipliers attached to the archimedian places of the field $K$ (see introduction for the case $K = \Q$, the appendix for the general case,  and \cite{W2},\cite{L1}).

 We wish to examine the following sum
 $$
 \res_{(0)}\xi_K(s)ds  +  \res_{(1)}\xi_K(s)ds - \lim_{A \to +\infty}\int\limits_{A - i\infty}^{A + i\infty}\xi_K(s)ds + \lim_{A \to -\infty}\int\limits_{A - i\infty}^{A + i\infty}\xi_K(s)ds
 $$
and to show that it is equal to zero\footnote{This question was discussed during my lectures in Moscow. The  ideas on how to deal with this sum belong to Volodja Zhgun and  Irina Rezvjakova. }.
This is actually the case (see appendix) and we see that our observation also has  some meaning for the number fields.
 We then have to describe  what stands for function spaces and maps from section 3 in the number field case.

Let $K/\Q$ be a finite extension of $\Q$ and $\A_K$  be the adelic group (= adelic product of all local fields $K_v$ for all places $v$).
We set
$$
K_v^* \supset {\hat{\cal O}}^*_v \quad\mbox{for finite (non-archimedean) places}~v,
$$
$$
K_v^* \supset {\hat{\cal O}}^*_v :  = \{ \pm 1  \}\quad\mbox{if}~K_v = \R,
$$
$$
K_v^* \supset {\hat{\cal O}}^*_v : = \{|z| = 1 \}\quad\mbox{if}~K_v = \C,
$$
and ${\hat{\cal O}}^* = \prod_v {\hat{\cal O}}^*_v$. This is the maximal compact subgroup in $\A^*_K$. We may take
$$
\Gamma_K = \oplus_v  K_v^*/{\hat{\cal O}}^*_v  =  \oplus_{\mbox{{\footnotesize finite}}~v} \Z \oplus_{\mbox{{\footnotesize infinite}}~v}\R = (\mbox{discrete group of ideals})\oplus (\mbox{vector space over}~\R).
$$
We now have to define  a space ${\cal S}(\Gamma_K)$ and  a map
$$
{\cal S}(\A_K)^{{\hat{\cal O}}^*} \stackrel{\sim}\rightarrow   {\cal S}(\Gamma_K).
$$
The space $ {\cal S}(\A_K)$ of Bruhat-Schwartz functions is well defined  in \cite{Br} (see also \cite{S, OP2}). On the local archimedean components $K_v$ of $\A_K$, these spaces are the Schwartz spaces and on the local non-archimedean components $K_v$ of $\A_K$, they are the Bruhat spaces ${\cal D}(K_v)$ (see the beginning of the introduction). For the local spaces ${\cal S}(\Gamma_v)$ one may choose the spaces
${\cal S}(K_v)^{{\hat{\cal O}}_v^*}$ themselves,
considering them  as  functional spaces  on $\Gamma_v$.

In the case of finite places, the groups  $ \Gamma_v$ are equal to $ \Z$ and we can develop the theory exactly as for the case of algebraic curves. Namely, we can introduce a holomorphic torus ${\mathbb T}_v$ and the same functional spaces ${\cal D}(\Gamma_v) \cong \C[{\mathbb T}_v]$ and ${\cal D}_+(\Gamma_v) \cong \C_+[{\mathbb T}_v]$ as above. The Fourier transform will be well defined in the same way but the numbers $q_v (=q_x)$ are  not the powers of a single  $q$. For example, if $K = \Q$ then the $v$ correspond to  all primes $p$ and $q_v = p$. For  this  reason and  (more importantly)  by the structure of  archimedian places we are forced to consider the function spaces on the groups  $\bar{\mathbb T}_v$, universal coverings of the tori ${\mathbb T}_v$. Certainly, ${\mathbb T}_v \cong \C$.
For simplicity, we now consider the case  $K = \Q$\footnote{The interested reader can easily extend this discussion  to a general number field using, for example, \cite{W2}[ch. VII, nn 2, 6], \cite{L1}[ch. VII, n  4].}.

Then for the space   ${\cal D}(\Q_p)^{\Z_p^*}$,  we have  the following space on $\C$:
$$
{\cal F}_p(\C) : = \{P(p^{-s})\vert P \in \mbox{the space of Laurent polynomials},~s \in \C  \}
$$
and, for the space ${\cal D}_+(\Q_p)^{{\cal O}_p^*}$,  the space
$$
{\cal F}_{p, +}(\C) : = \{ (1 - p^{-s})^{-1}P(p^{-s})\vert P \in  \mbox{the space of Laurent polynomials},~s \in \C  \}.
$$
A more complicated problem is what to choose for the infinite, archimedian places. In our case, $K_v = \R$. The largest possible choice is the space
${\cal S}(\R)$. The traditional choice will be the much smaller subspace
$$
{\cal D}_+(\R) = \{\mbox{linear span of}~x^n\exp(-ax^2),~x \in \R,~n \in \Z_{\ge 0},~a \in \R_+   \}
$$
and for the space of invariants under ${\hat{\cal O}}_v^* = \{\pm \}$ we restrict ourselves  to the even powers of  $x$'s\footnote{Sergey Gorchinskiy has asked what is an analogue of the space ${\cal D}(\Q_p)$? A good question!}.

The Fourier transform on the group $\R$ (see the introduction) will preserve the space  ${\cal D}_+(\R)$. Indeed, the function $\exp(-\pi x^2)$ is invariant
under the Fourier transform, so that  we get
$$
\Fo (\exp(-\pi ax^2 ))  = \frac{1}{\sqrt a}\exp(-\pi a^{-1}x^2).
$$
From the point of view of the analogy between the fields like $\Q_p,~{\mathbb F}_p((t))$ and $\R$, it is worth to compare this formula with the action of the Fourier
transform on the standard functions from the section 3.1: $\Fo(\delta_{\ge m}) = q^{-1/2k - m}\delta_{\ge -k-m}$ with $k = 0$.

More generally, for a polynomial $P(x)$ we have
$$
\Fo (P(x)\exp (-\pi ax^2)) = Q(x)\exp (-\pi a^{-1}x^2)
$$
with a new polynomial $Q(x)$. See the details in \cite{W2}[ch. VII, n 2, prop. 3].

Now, we have to use the classical Mellin transform ${\rm M}$:
$$
{\rm M}_v(f) = \int_{\R^*_+}\vert x\vert^s f(x) dx/\vert x \vert,\quad x \in K_v, \quad f \in {\cal D}_+(K_v).
$$
Using the standard integral presentation for the $\Gamma$-function, one gets  that, for the field $K_v = \R, v = \infty$,
$$
{\rm M}_{\infty}(x^{2n}\exp (-ax^2))  =  a^{-n}a^{-s/2}\Gamma(s/2 + n).
$$
This suggests the following definition
$$
{\cal F}_{\infty, +}(\C) : = \{\mbox{linear span of}~a^{-s/2}\Gamma(s/2 + n),\quad s \in \C, \quad n \in \Z_{\ge 0},\quad a \in \R_+ \}
\footnote{In classical analysis, the choice of smaller subspaces in ${\cal S}(\R)$, with stronger conditions for growth at infinity, implies
better analytical properties of the Mellin transform (the so called Paley-Wiener type theorems, see a clear short survey in \cite{Shilov}[ch. VII, n 3.4]).}.
$$

We  now proceed along the line developed for the case of algebraic curves. We  introduce  the whole space $\otimes_v{\cal F}_{v, +}(\C)$,
 then we map a single copy of $\C$ into the product of the $\C$'s over all places $v$ and, by restriction, we get a space
   ${\cal F}_{ +}(\C)$.
  Then, the functional equation for the zeta-function of the field $K$ \cite{L1},\cite{W2} should imply  the commutativity of the diagram
 \begin{equation}
\begin{CD}
{\cal D}_+(\A_K)^{{\hat{\cal O}}^*} @> \Fo >>  {\cal D}_+(\A_K)^{{\hat{\cal O}}^*}\\
@ V {\rm M}VV   @ V {\rm M}VV\\
{\cal F}_+(\C) @ >i^* >> {\cal F}_+(\C)
\end{CD}
\end{equation}
  where the global Mellin map ${\rm M}$ is a product of the local ones  ${\rm M}_v$ over all $v$ and $i: \C \rightarrow \C$ with $i(s) = 1 - s $.
One can  also show that the Poisson formula becomes the residue relation on
$s$-plane we have stated above. We want nevertheless to finish with the following

{\bf Question}. Does  there exist an analogue of the  group dual to the group ${\Gamma_K}$ for a  number field $K$ and how can one  realize the whole analysis working entirely on this group ?

%\subsection{Explicit formulas again}
\section{Appendix  (Irina Rezvjakova)
 Residues of the Dedekind zeta functions of number fields}

Let $\Q$ be the field of rational numbers,
$K \supset  \Q$ be its finite extension of degree  $d=r_1 + 2r_2$ and  $D $ be its  discriminant.

The Dedekind zeta-function  of the field  $K$ can be defined by the following formula for
${\rm Re}(s) >1$:
$$
\zeta_K (s) = \sum\limits_{\mathfrak{a}} \Nm(\mathfrak{a})^{-s} =
\sum\limits_{n=1}^{+\infty} \frac{a(n)}{n^s},
$$
where
$$
a(n)= \sum\limits_{\mathfrak{a}: \Nm(\mathfrak{a}) = n} 1.
$$
Here, the summation in the first equality is done over all nonzero integral  ideals  $\mathfrak{a}$
from  the ring of integers of the field $K$ and $\Nm(\mathfrak{a})$ denotes  norm of the ideal $\mathfrak{a}$.

Starting from  $\zeta_K (s)$ one can define  the $\xi$-function
$$
\xi_K (s) = G_1(s)^{r_1}G_2(s)^{r_2} \zeta_K (s),
$$
where
$$
G_1(s) = \pi^{-s/2}\Gamma(s/2), \quad G_2(s) =
(2\pi)^{1-s}\Gamma(s).
$$
(compare with the corresponding extension of Riemann's zeta-function in section 1).
This function is meromorphic on the whole complex plane with first order poles
at points  $s=0$ and $s=1$. It satisfies the following functional equation \cite{L1, W2}:
$$
\xi_K (s) = |D|^{\frac12 -s} \xi_K (1-s).
$$

Let $A>1$ and  $T$ be sufficiently large real numbers.
Consider the rectangle
$\Gamma = \gamma_1 \cup \gamma_2 \cup \gamma_3 \cup \gamma_4$:
$$
\begin{array}{llllllll}
\gamma_1 & = & \{A + it & : &  -T& \leq& t &\leq T \},\\
\gamma_2 & = & \{\sigma + iT& : & 1-A &\leq &\sigma &\leq A \},\\
\gamma_3 & = & \{1-A + it &: &  -T &\leq &t &\leq T \},\\
\gamma_4 & = & \{\sigma + iT& : &  1-A &\leq &\sigma &\leq A \},
\end{array}
$$
where $s=\sigma+it$.
\smallskip
\smallskip
By Cauchy's theorem we get
$$
\frac{1}{2\pi i} \int\limits_{\Gamma} \xi_K (s) ds = \res_{s=0}
\xi_K (s) + \res_{s=1} \xi_K (s).
$$
We want to show that for the fixed  $A$ and  $T\to +\infty$ the integrals over
the paths $\gamma_2, \gamma_4$ will tend to zero.
This can be done in four steps.

1) The Stirling formula for the Gamma function
implies that in the domain $|{\rm Re}(s)| \le a, |{\rm Im}(s)| \ge 1$
($s=\sigma +it$)
$$
|\Gamma(s)| \ll |t|^{\sigma-\frac12} e^{-\frac{\pi}{2}t},
$$
(see, for example, \cite{VK}[appendix, \S 3]).

2) In the right  half-plane, under  the condition $\mbox{Re}(s) \ge 1+\varepsilon$ (where
 $\varepsilon$ is an arbitrary small positive  number) there  holds
the estimate
$$
\vert \zeta_K (s)\vert  \le c = c(\varepsilon).
$$
Proof. Looking at  the Dirichlet series $\zeta_K (s) =  \sum_{\mathfrak{a}} \Nm(\mathfrak{a})^{-s}$ we
can rewrite it as
$$
\sum_{n \ge 1, n \in \Z} \sum_{\mathfrak{a}\cap \Z = (n)} \Nm(\mathfrak{a})^{-s}.
$$
This immediately reduces the problem to the corresponding estimate for the Riemann's zeta function where it is obvious.

3) From 1), 2) and the functional equation we get that  the inequality
$$
\vert \zeta_K (s)\vert  \le c_1 (\varepsilon) |t|^{d \left( \frac12 - \sigma
\right)}
$$
is valid  for $-a \le {\rm Re}(s) \le -\varepsilon$, $|{\rm Im}(s)| \ge
1$.

4)To get  a bound in  the critical strip one has to apply
the Phragmen-Lindel\"of principle.

To find the order of growth  of Dedekind's zeta-function in the domain
$-\varepsilon \le {\rm Re} s \le 1+\varepsilon$ one can use the  Phragmen-Lindel\"of principle:
{\it Let the function
$f$ be holomorphic in a neighborhood  of the strip $a\le {\rm Re} s \le b$ where the estimate $|f(s)| \ll e^{|s|^A}$
holds with some  $A \ge 0$. If
$$
|f(a+it)| \le c_a (1+|t|)^{\alpha},
$$
$$
|f(b+it)| \le c_b (1+|t|)^{\beta},
$$
then,  for  $a\le \sigma\le b$, we have
$$
|f(\sigma+it)| \le c_a^{l(\sigma)} c_b^{1-l(\sigma)} (1+|t|)^{\alpha
l(\sigma) +\beta (1-l(\sigma))},
$$
where  $l(\sigma)$~ is a linear function such that  $l(a)=1$, $l(b)=0$.}

We apply this statement to the holomorphic function  $f(s) = (s-1) \zeta_K
(s)$ in the strip with ends  $a= -\varepsilon$, $b = 1+\varepsilon$, with
parameters $\alpha = d(1/2 +\varepsilon)+1$, $\beta= 1$ and function
$l(\sigma) = -\frac{1}{1+2\varepsilon} \sigma +
\frac{1+\varepsilon}{1+2\varepsilon}$. Using the principle, we get that in the given strip
$$
|f(\sigma+it)| \ll_{\varepsilon} (1+|t|)^{1+ d(1/2 +\varepsilon)
\left( -\frac{1}{1+2\varepsilon} \sigma +
\frac{1+\varepsilon}{1+2\varepsilon}\right) } \ll_{\varepsilon}
(1+|t|)^{1+ d \left( - \frac{\sigma}{2} +
\frac{1+\varepsilon}{2}\right)}.
$$
Therefore, under $-\varepsilon \le \sigma \le 1+\varepsilon$, $|t|
\ge 1$ there is the estimate
$$
|\zeta_K (\sigma+it)| \ll_{\varepsilon} |t|^{ d \left(
\frac{1-\sigma}{2} + \frac{\varepsilon}{2}\right)}.
$$
Thus, in the domain ${\rm Re}(s) \ge -a$, $|t| \ge 1$
$$
|\zeta_K (\sigma+it)| \ll |t|^{ d \left(
\frac{1+a}{2}\right)  },
$$
5) Thus,  the integrals of the function $\xi_K (s)$ over the paths  $\gamma_2, \gamma_4$
can be bounded from above up to a constant depending only on  $A$:
$$
T^{A-\frac12} e^{-\frac{\pi}{2}T} T^{d(\frac{1}{2}+A)}  \ll
e^{-\frac{\pi}{2}T} T^{(d +1)(A+\frac{1}{2})} \to 0 \quad \text {when} \quad T\to
+\infty.
$$
Consequently, going to the limit where  $T\to +\infty$, we get the equality:
\begin{equation}\label{eq-res-sum}
\frac{1}{2\pi i} \int\limits_{(A)} \xi_K (s) ds - \frac{1}{2\pi i}
\int\limits_{(1-A)} \xi_K (s) ds= \res_{s=0} \xi_K (s) + \res_{s=1}
\xi_K (s),
\end{equation}
where $(A)$ denotes the line  ${\rm Re}(s)  = A$.

\bigskip

{\bf Example}. Let $K  = \Q (\sqrt{-D})$ be an imaginary quadratic field,
where $-D \le -3$.
Denote by  $h= h(-D)$ the class number of the field $K$ and by
$\omega$ the number of units (recall that $\omega =2$ when $D > 4$).
In this case, the  $\xi$-function is equal to
$$
\xi_K (s) = (2\pi)^{-s} \Gamma(s) \zeta_K (s)
$$
and satisfies
$$
\xi_K (s) =  D^{\frac12 -s} \xi_K(1-s).
$$

Let us find values of the residues:
$$
\res_{s=1} \xi_K (s) = \frac{1}{2\pi} \res_{s=1} \zeta_K (s) =
\frac{h}{\omega \sqrt{D}}.
$$
From the functional equation one obtains
$$
\res_{s=0} \xi_K (s) = \zeta_K (0) = \left( \frac{2\pi}{\sqrt{D}}
\right)^{-1} (-\res_{s=1} \zeta_K (s)) = -\frac{h}{\omega}.
$$
Now, we compute the values of integrals.  Applying the functional equation and then the Mellin transform,
we come to
\begin{align*}
&\frac{1}{2\pi i} \int\limits_{(1-A)} \xi_K (s) ds = \frac{1}{2\pi
i} \int\limits_{(1-A)} D^{\frac12 -s}\xi_K (1-s) ds = \frac{1}{2\pi
i} \int\limits_{(A)} D^{s-\frac12} \xi_K (s) ds \\
&= D^{-\frac12} \sum\limits_{n=1}^{+\infty} a(n)\frac{1}{2\pi i}
\int\limits_{(A)} \Gamma(s) \left( \frac{2\pi n}{D}\right)^{-s} ds =
D^{-\frac12} \sum\limits_{n=1}^{+\infty} a(n) e^{-\frac{2\pi n}{D}}
= D^{-\frac12} \sum\limits_{\mathfrak{a}} e^{-\frac{2\pi
\mbox{\footnotesize Nm}(\mathfrak{a})}{D}}.
\end{align*}
Also,
\begin{align*}
&\frac{1}{2\pi i} \int\limits_{(A)} \xi_K (s) ds =
\sum\limits_{n=1}^{+\infty} a(n) e^{-2\pi n} =
\sum\limits_{\mathfrak{a}} e^{-2\pi \mbox{\footnotesize Nm}(\mathfrak{a})}.
\end{align*}
We get that
$$
\frac{h}{\omega} + \sum\limits_{\mathfrak{a}} e^{-2\pi
\mbox{\footnotesize Nm}(\mathfrak{a})} = \frac{h}{\omega \sqrt{D}} + D^{-\frac12}
\sum\limits_{\mathfrak{a}} e^{-\frac{2\pi \mbox{\footnotesize Nm}(\mathfrak{a})}{D}}.
$$
This equality also follows  from  the functional equation for the corresponding theta-series.
Namely, if  ${\cal A}$ is a class of ideals, then for the function
$$
\theta_{{\cal A}} (z) = \frac{1}{\omega} +
\sum\limits_{\mathfrak{a} \in {\cal A}} e^{2\pi i z
\mbox{\footnotesize Nm}(\mathfrak{a})}, \quad \mbox{Im} z
>0,
$$
  the functional equation:
$$
\theta_{{\cal A}} (z) = \frac{i}{z\sqrt{D}}
\theta_{{\cal A}^{-1}} \left( -\frac{1}{Dz}\right)
$$
holds (see, for example, \cite{H2}). Thus, the function $\theta_K (z) = \sum\limits_{{\cal A}}
\theta_{{\cal A}} (z) = \frac{h}{\omega} +
\sum\limits_{\mathfrak{a}} e^{2\pi i z \mbox{\footnotesize Nm}(\mathfrak{a})}$ satisfies  the analogous
functional equation
$$
\theta_{K} (z) = \frac{i}{z\sqrt{D}} \theta_{K} \left(
-\frac{1}{Dz}\right).
$$
This equality taken  at  the point  $z=i$ will be exactly  the  equality which we  got from the Cauchy formula.

\end{document}